\documentclass[12pt]{article}

\usepackage{amssymb}

\usepackage{amsmath}

\usepackage[mathcal]{euscript}

\usepackage{amsthm}

\newcommand{\XX}{\mathbb{X}}
\newcommand{\YY}{\mathbb{Y}}
\newcommand{\e}{\varepsilon}

\begin{document}

\begin{center}
{\bf \large Asymptotic Expansions \\ for Stationary Distributions  of  \vspace{1mm} \\ Nonlinearly Perturbed  Semi-Markov Processes. II}
\end{center}
\vspace{2mm}

\begin{center}
{\large Dmitrii Silvestrov\footnote{Department of Mathematics, Stockholm University, SE-106 81 Stockholm, Sweden. \\ 
Email address: silvestrov@math.su.se} 
and Sergei Silvestrov\footnote{Division of Applied Mathematics, School of Education, Culture and Communication, M{\"a}lardalen  University, SE-721 23 V{\"a}ster{\aa}s, Sweden. \\ 
Email address: sergei.silvestrov@mdh.se}}
\end{center}
\vspace{2mm}

Abstract:
Asymptotic expansions with explicit upper bounds for remainders are given for stationary distributions  of nonlinearly perturbed semi-Markov processes with  finite phase spaces. The corresponding algorithms are based on a special technique of sequential phase space reduction, which can be applied to processes  with an arbitrary asymptotic communicative structure of phase spaces. \\
 


{\bf 1. Introduction}  \\

In this paper, we present new algorithms  for  construction of asymptotic expansions for stationary distributions of nonlinearly perturbed semi-Markov processes with a finite phase space.  This is Part II of the paper,  Silvestrov, D. and Silvestrov,  S. (2016), where algorithms  for  constructing of  asymptotic expansions with remainders of a standard form $o(\cdot)$ have been given. In Part II, we present algorithms  for  construction  asymptotic expansions of a more advanced form, with explicit upper bounds for remainders.

We consider models, where  the phase space is one class of communicative states, for embedded Markov chains of pre-limiting perturbed semi-Markov processes, while it can possess an arbitrary communicative structure, i.e., can consist of one or several closed classes of communicative states and, possibly,  a class of transient states,  for the limiting embedded Markov chain. 

The initial perturbation conditions are formulated  in the forms of Taylor and Laurent asymptotic expansions with explicit upper bounds for remainders, respectively,  for transition probabilities (of  embedded Markov chains)  and  expectations of sojourn times, for perturbed semi-Markov processes.  

The algorithms are based  on special time-space screening procedures for sequential phase space reduction and algorithms for  re-calculation of asymptotic expansions with explicit upper bounds for remainders, which constitute perturbation conditions for the  semi-Markov processes with reduced phase spaces. 

The final asymptotic expansions for stationary distributions of nonlinearly perturbed semi-Markov processes are given in the form of Taylor asymptotic expansions, with explicit upper bounds for remainders.  

The algorithms presented in the paper have an universal character. They can be applied to perturbed semi-Markov processes with an arbitrary asymptotic communicative structure and are computationally effective due to recurrent character of computational procedures.

The survey of works in the area and detailed comments concerned with the proposed method are given in the Introduction to Part I of this paper. 

Here, we would like, only, to mention books, which contains parts devoted to perturbed Markov chains and  semi-Markov processes,  and problems related to asymptotic expansions for the above  models. These are,  Korolyuk and Turbin  (1976, 1978), Courtois  (1977),  Seneta (1981, 2006),   Stewart and  Sun (1990),  Kartashov (1996),  Stewart (1998, 2001), Yin and  Zhang (1998, 2005, 2013), Korolyuk, V.S. and Korolyuk, V.V. (1999),   Konstantinov,   Gu,  Mehrmann and Petkov (2003), Bini, Latouche and Meini (2005), Korolyuk and  Limnios (2005),   Gyllenberg and Silvestrov (2008), and  Avrachenkov,  Filar and Howlett  (2013). 

A comprehensive bibliography of works in the area can be found in these books and, also, in the research report by Silvestrov, D. and Silvestrov, S. (2015), which is an extended preliminary version of the present paper.

In conclusion, we would like to mention that, by our opinion, the results presented in the paper have a good potential for continuation of studies. We comment some prospective directions for future studies in the last section of the paper. 

Part II includes four sections and two appendices. In Section 2, we present so-called operational rules for  Laurent asymptotic expansions with explicit upper bounds for remainders. In Section 3, we present basic perturbation conditions and algorithms  for  construction of  asymptotic expansions with explicit upper bounds for remainders, for transition characteristics of nonlinearly perturbed semi-Markov processes with reduced phase spaces. In Section 4, we present algorithms for construction of asymptotic expansions with explicit upper bounds for remainders, for expected hitting times for nonlinearly perturbed semi-Markov processes. In Section 5, we   present an algorithm for construction of asymptotic expansions with explicit upper bounds for remainders, for stationary distributions of nonlinearly perturbed semi-Markov processes. In Appendix A, we give proofs of lemmas representing operational rules for  Laurent asymptotic expansions without and with explicit upper bounds for remainders. In Appendix B, we discuss and  present  examples  illustrating algorithms for construction of asymptotic expansions for stationary distributions of nonlinearly perturbed semi-Markov processes developed in the present paper. \\

{\bf 2. Laurent asymptotic expansions  with explicit  upper bounds \\ \makebox[10mm]{} for remainders} \\

In this section, we present so-called operational rules for  Laurent asymptotic expansions with explicit upper bounds for remainders. The corresponding proofs and comments are given in Appendix A.

Let $A(\e)$ be a real-valued function 
defined on an interval $(0, \e_0]$, for some $0 < \e_0 \leq 1$, and given on this interval
by a  Laurent asymptotic expansion,
\begin{align}\label{expap}
A(\e) = a_{h_A}\e^{h_A} + \cdots + a_{k_A}\e^{k_A} + o_A(\e^{k_A}),
\end{align}
where {\bf (a)} $- \infty < h_A \leq k_A < \infty$ are integers, {\bf
(b)} coefficients $a_{h_A}, \ldots, a_{k_A}$ are real numbers, {\bf (c)} 
$|o_A(\e^{k_A})|$ $\leq G_A \e^{k_A + \delta_A}$,  for $0 < \e \leq \e_A$, where {\bf (d)}  $0 < \delta_A \leq 1, 0 < G_A < \infty$ and  $0 < \e_A \leq \e_0$.
 
We refer to such  Laurent asymptotic expansion as a $(h_A, k_A, \delta_A, G_A,  \e_A)$-expansion.

The   $(h_A, k_A, \delta_A, G_A,  \e_A)$-expansion is also a $(h_A, k_A)$-expansion,  according the definition given in Part I of the paper, since,  $o_A(\e^{k_A})/\e^{k_A} \to 0$ as $\e \to 0$.  

We say that  $(h_A, k_A, \delta_A, G_A,  \e_A)$-expansion $A(\e)$ is pivotal if it is known that \mbox{$a_{h_A} \neq 0$.}

It is useful to note that there is no sense to consider, it seems,  a more general case of upper bounds for the remainder $o_A(\e^{k_A})$, with parameter $\delta_A > 1$. Indeed, let us define $k'_A =  k_A +  [\delta_A] - {\rm I}(\delta_A = [\delta_A])$ and $\delta'_A = \delta_A -  [\delta_A] + {\rm I}(\delta_A = [\delta_A])
\in (0, 1]$. 
The $(h_A, k_A, \delta_A, G_A,  \e_A)$-expansion $A(\e)$ can be re-written in the equivalent form of the  
$(h_A, k'_A, \delta'_A, G_A, \e_A)$-expansion,
$A(\e) = a_{h_A}\e^{h_A} + \cdots + a_{k_A}\e^{k_A} + 0 \e^{k_A +1} + \cdots + 0 \e^{k'_A} +  o'_A(\e^{k'_A})$,
with the remainder  $o'_A(\e^{k'_A}) = o_A(\e^{k_A})$,  which satisfies inequalities $|o'_A(\e^{k'_A})| = |o_A(\e^{k_A})| \leq 
 G_A \e^{k_A + \delta_A} = G_A \e^{k'_A + \delta'_A}$, for $0 < \e \leq \e_A$.
 
The above remarks imply that the  asymptotic expansion $A(\e)$ can be represented in different forms. In such cases, we consider  forms with larger parameters $h_A$ and $k_A$ as more informative. As far as parameters  $\delta_A, G_A$ and $\e_A$ are concerned, we consider as more informative forms, first, with larger values of parameter $\delta_A$, second, with smaller values of parameter $G_A$  and, third, with larger values of parameter $\e_A$.   

In what follows, lemmas, theorems and relations from  Part I of the paper are indexed  by symbol $*$.

The following proposition supplements Lemma 1$_*$.  \vspace{1mm}

{\bf Lemma 1}. {\em 
If $A(\e) = a'_{h'_A}\e^{h'_A} + \cdots + a'_{k'_A} \e^{k'_A} + o'_A(\e^{k'_A}) =  a''_{h''_A}\e^{h''_A} + \cdots + a''_{k''_A} \e^{k''_A} + o''_A(\e^{k''_A}), \e \in (0, \e_0]$ can be represented as, respectively,  $(h'_A, k'_A,  \delta'_A, G'_A$, $\e'_A)$- and $(h''_A$, $k''_A, \delta''_A, G''_A, \e''_A)$-expansion, then 
the $(h_A, k_A)$-expansion $A(\e) = a_{h_A}\e^{h_A} + \cdots + a_{k_A}\e^k + o_A(\e^{k_A}), \e \in (0, \e_0]$  given in Lemma 1$_*$ is an  $(h_A, k_A, \delta_A$, $G_A, \e_A)$-expansion,  with parameters $G_A, \delta_A$ and $\e_A$ chosen in the following way consistent with the priority order 
described above{\rm :}
\begin{equation*}\label{equar}
(\delta_A, G_A, \e_A) = \left\{
\begin{array}{ll}
(\delta''_A, G''_A,  \e''_A)  & \ \text{if} \   k'_A < k''_A \, or \, k'_A = k''_A, \\
&  \quad \ \delta'_A < \delta''_A,  \\
(\delta'_A = \delta''_A, G'_A \wedge G''_A,  \e'_A \wedge \e''_A) & \ \text{if} \   k'_A = k''_A,  \delta'_A = \delta''_A,  \\
(\delta'_A,  G'_A,  \e'_A)   & \ \text{if} \   k' _A > k''_A \, or \, k'_A = k''_A, \\
&  \quad \ \delta'_A > \delta''_A.
\end{array}
\right.
\end{equation*} 
} \vspace{1mm}

It is also useful to mention that  a constant $a$ can be interpreted  as function $A(\e) \equiv a$. Thus, $0$ can be represented, for any integer $- \infty < h \leq k < \infty$, as the $(h, k, \delta_{h, k}, G_{h, k}$, $\e_{h, k})$-expansion, $0 = 0 \e^h + \ldots + 0 \e^k + o(\e^k)$, with  remainder $o(\e^k) \equiv 0$ and, thus, parameters $\delta_{h, k} = 1$, an arbitrary small $G_{h, k} \in (0, \infty)$, and $\e_{h, k} = \e_0$. Also,  $1$ can be represented, for any integer $0 \leq k < \infty$, as the $(0, k, \delta_k, G_k, \e_k)$-expansion,  $1 = 1 + 0 \e + \ldots + 0 \e^k + o(\e^k)$, with remainder $o(\e^k) \equiv 0$ and, thus, parameters $\delta_k = 1$, an arbitrary small  $G_k \in (0, \infty)$,  and $\e_k = \e_0$.

Let us consider four Laurent asymptotic  expansions, $A(\e) = a_{h_A}\e^{h_A} +
\cdots + a_{k_A}\e^{k_A} + o_A(\e^{k_A})$, $B(\e) = b_{h_B}\e^{h_B} +
\cdots + b_{k_B}\e^{k_B} + o_B(\e^{k_B})$, $C(\e) = c_{h_C}\e^{h_C} +
\cdots + c_{k_C}\e^{k_C} + o_C(\e^{k_C})$, and $D(\e) =
d_{h_D}\e^{h_D} + \cdots + d_{k_D}\e^{k_D} + o_D(\e^{k_D})$ defined for $0 < \e \leq \e_0$, for some $0 < \e_0 \leq 1$.

The following Lemma presents   operational rules for computing parameters of upper bounds for remainders of  Laurent asymptotic expansions.  \vspace{1mm}

{\bf Lemma 2}. {\em 
The above asymptotic expansions
have the following operational rules for computing remainders{\rm :}

{\bf (i)} If $A(\e), \e \in (0, \e_0]$ is a $(h_A, k_A,  \delta_A,  G_A,\e_A)$-expansion
  and $c$ is a constant, then $C(\e) = cA(\e), \e \in (0, \e_0]$ is a $(h_C, k_C, \delta_C, G_C,  \e_C)$-expansion with parameters $h_C, k_C$ and coefficients $c_r, r = h_C, \ldots, k_C$  given in proposition  {\bf (i)} of Lemma 2$_*$, and parameters $\delta_C, G_C$ and  $\e_C$ given by the formulas{\rm :}  

{\bf (a)} $\delta_C = \delta_A${\rm ;}  

{\bf (b)}  $G_C = |c|G_A${\rm ;}  

{\bf (c)} $\e_C =  \e_A$.

{\bf (ii)} If $A(\e), \e \in (0, \e_0]$ is a $(h_A, k_A, \delta_A,  G_A, \e_A)$-expansion
  and $B(\e), \e \in (0, \e_0]$ is a $(h_B, k_B$, $\delta_B, G_B,  \e_B)$-expansion, then $C(\e) =
  A(\e)+ B(\e)$, $\e \in (0, \e_0]$ is a $(h_C, k_C,  \delta_C$, $G_C, \e_C)$-expansion with parameters $h_C,  k_C$ and coefficients $c_{r}, r = h_C, \ldots, k_C$ given in proposition  {\bf (ii)} of Lemma 2$_*$, and parameters $\delta_C, G_C$ and   $\e_C$ given by formulas{\rm :} 
  
 {\bf (a)} $\delta_C = \delta_A {\rm I}(k_A < k_B) \  + \
(\delta_A \wedge \delta_B) {\rm I}(k_A = k_B) \ + \ \delta_B {\rm I}(k_B <  k_A) \geq \delta_A \wedge \delta_B${\rm ;}   

{\bf  (b)} $G_C =  G_A \e_C^{k_A \, +  \, \delta_A - k_C - \delta_C}  \, + \, \sum_{k_C < i \leq k_A} |a_i|\e_C^{i - k_C - \delta_C}  
+  \, G_B \e_C^{k_B + \delta_B - k_C - \delta_C}  \\
\makebox[10mm]{}   +  \, \sum_{k_C < j \leq k_B} |b_j|\e_C^{j - k_C - \delta_C}${\rm ;}   

{\bf  (c)} $\e_C = \e_A \wedge \e_B$.
 
{\bf (iii)} If $A(\e), \e \in (0, \e_0]$ is a $(h_A, k_A, \delta_A, G_A, \e_A)$-expansion
  and $B(\e), \e \in (0, \e_0]$ is a $(h_B, k_B$, $\delta_B, G_B, \e_B)$-expansion, 
  then $C(\e) = A(\e) \cdot B(\e)$,  $\e \in (0, \e_0]$ is a $(h_C, k_C,  \delta_C$, $G_C, \e_C)$-expansion with parameters $h_C, k_C$ and coefficients $c_{r}, r = h_C, \ldots, k_C$ given in proposition  {\bf (iii)} of Lemma 2$_*$, and parameters $\delta_C, G_C$ and $\e_C$ given by formulas{\rm :} 
  
 {\bf  (a)} $\delta_C = \delta_A {\rm I}(k_A + h_B < k_B + h_A ) \ + \ (\delta_A \wedge \delta_B) {\rm I}(k_A + h_B = k_B + h_A) \\ 
 \makebox[10mm]{}  + \ \delta_B {\rm I}(k_A + h_B > k_B + h_A) \geq \delta_A \wedge \delta_B${\rm ;}   
  
 {\bf (b)}  $G_C =  \sum_{k_C < i + j, h_A \leq i \leq k_A, h_B \leq j \leq  k_B} |a_i| |b_j| \e_C^{i + j -  k_C - \delta_C} \\ 
 \makebox[10mm]{}   + \ G_A \sum_{h_B \leq j \leq k_B} |b_j| \e_C^{j + k_A + \delta_A  - k_C  - \delta_C}   
\ + \ G_B  \sum_{h_A \leq i \leq k_A} |a_i|  \e_C^{i + k_B  + \delta_B - k_C  - \delta_C} \\
\makebox[10mm]{}  + \ G_A G_B  \e_C^{k_A+ k_B  + \delta_A  +\delta_B - k_C  - \delta_C}${\rm ;} 

{\bf  (c)} $\e_C = \e_A \wedge \e_B$. 

{\bf (iv)} If $B(\e), \e \in (0, \e_0]$ is a pivotal $(h_B, k_B, \delta_B, G_B, \e_B)$-expansion, then there exists 
$\e_C \leq  \e'_0 \leq \e_0$ such that 
 $B(\e) \neq 0, \e \in (0, \e'_0]$, and $C(\e) = \frac{1}{B(\e)}, \e \in (0, \e'_0]$ is a pivotal $(h_C, k_C, \delta_C, G_C, \e_C)$-expansion with parameters $h_C, k_C$ and coefficients $c_{r}, r = h_C, \ldots, k_C$ given in proposition  {\bf (iv)} of 
Lemma 2$_*$, and parameters $\delta_C, G_C$ and  $\e_C$ given formulas{\rm :}

{\bf  (a)} $\delta_C = \delta_B${\rm ;} 

{\bf (b)} $G_C =  \ (\frac{|b_{h_B}|}{2})^{-1} \big( \sum_{k_B - h_B < i + j, h_B \leq i \leq k_B, h_C \leq j \leq  k_C} |b_i| |c_j| \e_C^{i + j - k_B + h_B - \delta_B} \\\makebox[10mm]{}  + \ G_B \sum_{h_C \leq j \leq k_C} |c_j| \e_C^{j + h_B} \big)${\rm ;} 

{\bf  (c)} $\e_C = \e_B \wedge \tilde{\e}_B$, where $\tilde{\e}_B =  
\big( \frac{|b_{h_B}|}{2( \sum_{h_B < i \leq k_B} |b_i| \e_B^{i - h_B - \delta_B} )
+  \ G_B \e_B^{k_B  - h_B})} \big)^{\frac{1}{\delta_B}}$.

{\bf (v)} If $A(\e), \e \in (0, \e_0]$ is a $(h_A, k_A,  \delta_A, G_A, \e_A)$-expansion, 
$B(\e), \e \in (0, \e_0]$ is a pivotal $(h_B, k_B, \delta_B, G_B, \e_B)$-expansion, then there exists 
$\e_D \leq  \e'_0 \leq \e_0$ such that  $B(\e) \neq 0, \e \in (0, \e'_0]$, and  $D(\e) = 
\frac{A(\e)}{B(\e)}$ is a $(h_D, k_D$,  $\delta_D, G_D$, $\e_D)$-expansion with parameters $h_D, 
k_D$ and  coefficients $d_{r}, r = h_D, \ldots, k_D$ 
given in proposition  {\bf (v)} of Lemma 2$_*$, 
and parameters $\delta_D, G_D, \e_D$ given by formulas{\rm :} 

{\bf  (a)} $\delta_D = \delta_A {\rm I}( h_C+ k_A < h_A + k_C) \ + \ (\delta_A \wedge \delta_C)  {\rm I}( h_C+ k_A = h_A + k_C) \\ 
\makebox[10mm]{}  + \ \delta_C {\rm I}( h_A+ k_C < h_C + k_A) \geq \delta_A \wedge \delta_B${\rm ;}   

{\bf  (b)} $G_D =   \sum_{k_D < i + j, h_A \leq i \leq k_A, h_C \leq j \leq  k_C} |a_i| |c_j| \e_D^{i + j -  k_D - \delta_D} \\
\makebox[10mm]{}  + \ G_A \sum_{h_C \leq j \leq k_C} |c_j| \e_D^{j + k_A + \delta_A  - k_D  - \delta_D}  
\ + \ G_C  \sum_{h_A \leq i \leq k_A} |a_i|  \e_D^{i + k_C  + \delta_C - k_D  - \delta_D}  \\
\makebox[10mm]{}  + \  G_A  G_C \e_D^{k_A + k_C  + \delta_A + \delta_C - k_D  - \delta_D}${\rm ;}   

{\bf (c)} $\e_D = \e_A \wedge \e_C$, \\
where coefficients $c_{r}, r = h_C, \ldots, k_C$ and parameters $h_C, k_C, \delta_C, G_C, \e_C$ are given for the  $(h_C, k_C, \delta_C, G_C, \e_C)$-expansion  of function $C(\e) = \frac{1}{B(\e)}$  in the above proposition {\bf (iv)}, or by formulas{\rm :} 

{\bf  (d)} $\delta_D = 
\delta_A {\rm I}(k_A - h_B < k_B  - 2 h_B + h_A) \ + \
(\delta_A \wedge \delta_B) {\rm I}(k_A - h_B = \\  \makebox[12mm]{}   k_B  - 2 h_B + h_A)  
+ \ \delta_B {\rm I}(k_A - h_B > k_B  - 2 h_B + h_A) \geq \delta_A \wedge \delta_B${\rm ;}    

{\bf  (e)}  $G_D =  (\frac{|b_{h_B}|}{2})^{-1} \big( \sum_{k_A \wedge (h_A + k_B - h_B) < i + j, h_A \leq i \leq k_A, h_D \leq j \leq  k_D} |a_i| \\
\makebox[10mm]{}   \times |d_j|  \e_D^{i + j  - k_D - h_B - \delta_D}  + \ \sum_{k_A \wedge (h_A + k_B - h_B) < i \leq k_A} |a_i| \e_D^{i  - h_B - k_D  - \delta_D} \ \\  
\makebox[10mm]{}  + \  G_A \e_D^{k_A \delta_A  - h_B - k_D - \delta_D}  +  \ G_B \sum_{h_D \leq j \leq k_D} |d_j|\e_D^{j + k_B + \delta_B - h_B  - k_D - \delta_D}  \big)${\rm ;}  

{\bf  (f)} $\e_D = \e_A \wedge  \e_B \wedge \tilde{\e}_B$, where $\tilde{\e}_B$ is given in relations {\bf (c)} of proposition {\bf (iv)}. 
} \vspace{1mm} 

{\bf Remark 1}. Coefficients $\e_B, \e_C, \e_D \in (0, 1]$ are taken  to nonnegative powers in all terms penetrating the sums,  which define parameters $G_C$, $G_D$ and $\tilde{\e}_B$,  in Lemma 2.  

The following operational rules for computing remainders for multiple summation and multiplication of Laurent asymptotic  expansions, used in what follows, are analogues of the corresponding summation and multiplication rules given in Lemma 2. \vspace{1mm}

{\bf Lemma 3}. {\em 
Let  $A_m(\e) =  a_{h_{A_{m}}, m}\e^{h_{A_{m}}} +
\cdots + a_{k_{A_{m}}, m}\e^{k_{A_{m}}} + o (\e^{k_{A_{m}}}), \e \in (0, \e_0]$ be a $(h_{A_m},  k_{A_m}, \delta_{A_m}, G_{A_m}, \e_{A_m})$-expansion, for  $m = 1, \ldots, N$. In 
this case{\rm :} 

{\rm \bf (i)} $B_n(\e) = A_1(\e) + \cdots + A_n(\e), \e \in (0, \e_0]$ is, for every $n = 1, \ldots, N$,   a $(h_{B_n}, k_{B_n}$, $\delta_{B_n}, G_{b_n}, \e_{B_n})$-expansion, with parameters $h_{B_n}, \, k_{B_n}, \, n = 1, \ldots, N$ and coefficients $b_{h_{B_n} + l, n}, l = 0, \ldots,   k_{B_n} - h_{B_n}, n = 1, \ldots, N$ given in proposition {\bf (i)} of Lemma 3$_*$, and parameters $G_{B_n}, \delta_{B_n}$, $\e_{B_n}, n = 1, \ldots, N$ given by formulas{\rm :}

{\bf (a)} $\delta_{B_n} =    \min_{m \in {\mathbb K}_n} \delta_{A_m} \geq \delta^*_N = \min_{ 1 \leq m \leq n} \delta_{A_m}$,  
where ${\mathbb K}_n = \{m:  1 \leq m \leq n,$ \\ 
\makebox[10mm]{}  $k_m = \min(k_1, \ldots, k_n) \}${\rm ;} 

{\bf (b)} $G_{B_n} =   \sum_{1 \leq i \leq n} \big( G_{A_i} \e_{B_{n}}^{k_{A_i} + \delta_{A_i} - k_{B_{n}}  - \delta_{B_{n}} }  
+  \sum_{k_{B_{n}}  < j \leq k_{A_i}} |a_{A_{i},j}| \e_{B_{n}}^{j - k_{B_{n}} - \delta_{B_{n}}} \big)${\rm ;} 

{\bf (c)} $\e_{B_n} = \min(\e_{A_1}, \ldots, \e_{A_n})$. 

{\bf (ii)} $C_n(\e) = A_{1}(\e) \times \cdots \times  A_n(\e), \, \e \in (0, \e_0]$ is, for $n = 1, \ldots, N$,   a $(h_{C_n}, k_{C_n},  \delta_{C_n}$, 
$G_{C_n}, \e_{C_N})$-expansion
with parameters $h_{C_n}, k_{C_n}, n = 1, \ldots, N$ and coefficients $c_{h_{C_n} + l, n}$, $l = 0, \ldots,   k_{C_n} - h_{C_n}, n = 1, \ldots, N$ given in proposition {\bf (ii)} of Lemma 3$_*$, and parameters $G_{C_n}, \delta_{C_n}$, $\e_{C_n}, n = 1, \ldots, N$ given by formulas{\rm :}

{\bf (a)} $\delta_{C_n} = \min_{m \in {\mathbb L}_n} \delta_{A_m} \geq  \delta^*_N$, 
where ${\mathbb L}_n  = \{m:  1 \leq m \leq n, \,  (k_{A_m}  \\   
\makebox[10mm]{} + \sum_{1 \leq r \leq n, r \neq m} h_{A_r}) = \min_{1 \leq l \leq n} (k_{A_l} + \sum_{1 \leq r \leq n, r \neq l} h_{A_r}) \}${\rm ;} 

{\bf (b)} $G_{C_n} \ =  \ \sum_{k_{C_n} < l_1 + \cdots + l_n, h_{A_i} \leq l_i \leq k_{A_{i}}, 1 \leq i \leq n} 
\ \prod_{1 \leq i \leq n} \, |a_{A_{i}, l_i}| \e_{C_n}^{l_1 + \cdots + l_n -  k_{C_n} - \delta_{C_n}} \\  
\makebox[10mm]{} + \ \sum_{1 \leq j \leq n} \prod_{1 \leq i \leq n, i \neq j}  \big( \sum_{h_{A_i} \leq l \leq k_{A_i}} |a_{A_{i}, l}| \e_{C_{n}}^{l}  \\
\makebox[10mm]{} + \ G_{A_i} \e_{C_{n}}^{k_{A_{i}}
 +  \delta_{A_{i}}} \big) G_{A_{j}}  \e_{C_{n}}^{k_{A_j} + \delta_{A_j} - k_{C_n}  - \delta_{C_n}}${\rm ;} 
 
{\bf (c)} $\e_{C_n} = \min_{1 \leq i \leq n} \e_{A_{i}}$. 

{\bf (iii)} Parameters $\delta_{C_n}, G_{C_n}, \e_{C_n}, n = 1, \ldots, N$ in upper bounds for remainders in the asymptotic expansions for functions $B_n(\e) = A_1(\e) + \cdots + A_n(\e), n = 1, \ldots, N$ and $C_n(\e) = A_1(\e) \times \cdots \times A_n(\e), n = 1, \ldots, N$ are invariant with respect to any permutation, respectively, of summation and multiplication order in the above formulas. } \vspace{1mm}

The summation and multiplication rules for computing of upper bounds for remainders given in  propositions  {\bf (ii)} and {\bf (iii)} of Lemma  2 possess the  communicative property, but do not possess the associative and distributional properties. 

Lemma 2 let us get an effective low bound for parameter $\delta_A$ for any $(h_{A}, k_{A}, \delta_{A}$, $G_{A}, \e_{A})$-expansion $A(\e)$ obtained as the result of a  finite sequence of operations (described in Lemma 2) performed over expansions from some finite set of such expansions.

The following lemma summarize these properties of Laurent asymptotic expansions with explicit upper bounds for remainders. \vspace{1mm}

{\bf Lemma 4}. {\em  
The summation and multiplication operations for Laurent asymptotic expansions defined in Lemma 2 possess the following algebraic properties, which should be understood as equalities for the corresponding  parameters of upper bounds for their remainders{\rm :} 

{\bf (i)} The functional identity, $C(\e) \equiv A(\e) +  B(\e) \equiv B(\e) +  A(\e)$, implies that $\delta_C = \delta_{A + B} = \delta_{B + A}, G_C = G_{A + B} = 
G_{B + A}$ and  $\e_C = \e_{A + B}  = \e_{B + A}$. 

{\bf (ii)} The functional identity,  $C(\e) \equiv A(\e) \cdot  B(\e) \equiv B(\e) \cdot  A(\e)$, implies that   $\delta_C = \delta_{A \cdot B} = \delta_{B \cdot A}, G_C = G_{A \cdot B} = 
G_{B \cdot A}$  and $\e_C = \e_{A \cdot B} = \e_{B \cdot A}$.

{\bf (iii)}  If  $A(\e)$  is $(h_{A}, k_{A}, \delta_{A}, G_{A}, \e_{A})$-expansion obtained as the result of a finite sequence of operations {\rm (}multiplication by a constant,  summation, multiplication, and division{\rm )} performed over  $(h_{A_i}, k_{A_i}, \delta_{A_i}, G_{A_i}, \e_{A_i})$-expansions  $A_i(\e), i = 1, \ldots, N$,  according the rules presented in Lemmas 2$_*$ and 2,  then $\delta_A \geq \delta^*_N = \min_{1 \leq i \leq N} \delta_{A_i}$. This  makes it possible to rewrite  $A(\e)$ as the $(h_{A}, k_{A}, \delta^*_N, G^*_{A, N}, \e_{A})$-expansion, with parameter $G^*_{A, N} = G_{A} \e_A^{\delta_A - \delta^*_N}$.
} \\

{\bf 3. Asymptotic expansions for transition characteristics of non- \\
\makebox[10mm]{} linearly perturbed semi-Markov processes with  reduced phase \\ 
\makebox[10mm]{}  spaces} \\

Let us recall the perturbed semi-Markov processes $\eta^{(\e)}(t), t \geq 0$, with phase space 
$\XX = \{1, \ldots, N \}$ and transition probabilities $Q^{(\e)}_{ij}(t), t \geq 0,\,  i, j \in \XX$,  introduced 
in Part I of the paper. These processes depend on a perturbation parameter $\e \in (0, \e_0]$, for some 
$0 < \e_0 \leq 1$.  We also recall transition probabilities of the corresponding embedded Markov chains, $p_{ij}(\e) = Q^{(\e)}_{ij}(\infty), i, j \in \XX$, and expectations of sojourn times  $e_{ij}(\e) = \int_0^\infty t  Q^{(\e)}_{ij}(dt)$, $i, j \in \XX$. 

We assume that condition ${\bf A}$,  introduced in Part I,  holds  for semi-Markov processes  $\eta^{(\e)}(t)$. In particular, we recall the transition sets $\YY_i, \, i \in \XX$ (which include states $j \in \XX$ with non-zero probabilities $p_{ij}(\e)$ and guarantee ergodicity of the processes $\eta^{(\e)}(t)$) introduced in this condition.

However, we replace  the perturbation condition ${\bf D}$  by the following stronger condition, in which the corresponding Taylor asymptotic expansions are given in the form with explicit upper bounds for remainders:
\begin{itemize}
\item [${\bf D'}$:] $p_{ij}(\e) =  \sum_{l = l_{ij}^-}^{l_{ij}^+} a_{ij}[l]\e^l + o_{ij}(\e^{l_{ij}^+}), \e \in (0, \e_0]$, for  $j \in \YY_i, i \in \XX$, where {\bf (a)}  $a_{ij}[ l_{ij}^-] > 0$ and $0 \leq l_{ij}^- \leq l_{ij}^+ < \infty$, for $j \in \YY_i, i \in \XX$; {\bf (b)}  $|o_{ij}(\e^{l_{ij}^+})| \leq G_{ij}\e^{l_{ij}^+ + \delta_{ij}}, 0 < \e \leq \e_{ij}$, for $j \in \YY_i, i \in \XX$, where $ 0 < \delta_{ij} \leq 1, 0 < G_{ij} < \infty$ and  $0 < \e_{ij} \leq \e_0$. 
\end{itemize}

Also, we replace  the perturbation condition ${\bf E}$  by the following stronger condition, in which the corresponding Laurent asymptotic expansions are given in the form with explicit upper bounds for remainders:
\begin{itemize}
\item [${\bf E'}$:] $e_{ij}(\e) =  \sum_{l = m_{ij}^-}^{m_{ij}^+} b_{ij}[l]\e^l + \dot{o}_{ij}(\e^{m_{ij}^+}), \e \in (0, \e_0]$, for $j \in \YY_i, i \in \XX$, where {\bf (a)} $b_{ij}[ l_{ij}^-] > 0$ and $- \infty <  m_{ij}^- \leq m_{ij}^+ < \infty$, for $j \in \YY_i, i \in \XX$; {\bf (b)}  $|\dot{o}_{ij}(\e^{l_{ij}^+})| \leq \dot{G}_{ij}\e^{m_{ij}^+ + \dot{\delta}_{ij}}, 0 < \e \leq \dot{\e}_{ij}$, for $j \in \YY_i, i \in \XX$, where $ 0 < \dot{\delta}_{ij} \leq 1, 0 < \dot{G}_{ij} < \infty$ and  $0 < \dot{\e}_{ij} \leq \e_0$. 
\end{itemize} 

As was pointed out in Part I, condition ${\bf A}$ implies that sets $\YY^+_{rr} = \YY_{rr} \setminus \{ r \}   \neq  \emptyset, r \in \XX$ and the non-absorption probability $\bar{p}_{rr}(\e) = 1 - p_{rr}(\e)  \in (0, 1]$, for $r \in \XX, \, \e \in (0, \e_0]$. This probability satisfy 
the following relation, for every $r \in \XX, \, \e \in (0, \e_0]$,
\begin{equation}\label{hopte}
\bar{p}_{rr}(\e) = 1 - p_{rr}(\e) = \sum_{j \in \YY^+_{rr}} p_{rj}(\e).
\end{equation}
 
The above relation let us construct an algorithm for getting  asymptotic expansions with explicit upper bounds for remainders, for non-absorption probabilities $\bar{p}_{rr}(\e)$.  \vspace{1mm} 

{\bf Lemma 5}. {\em 
Let conditions ${\bf A}$ and ${\bf D'}$  hold. Then, for every $r \in \XX$, the pivotal $(\bar{l}^-_{rr}, \bar{l}^+_{rr})$-expansion for the 
non-absorption probability $\bar{p}_{rr}(\e)$ given in Lemma 8$_*$ is, also, a $(\bar{l}^-_{rr}, \bar{l}^+_{rr},  \bar{\delta}_{rr}, \bar{G}_{rr}, \bar{\e}_{rr})$-expansion, with parameters $\bar{\delta}_{rr}, \bar{G}_{rr}$ and 
$\bar{\e}_{rr}$, which can be computed according the algorithm  described below,  in the proof of the lemma. 
} \vspace{1mm} 

{\bf Proof}. Let $r \in \YY_r$. First, propositions {\bf (i)} of Lemmas 3$_*$ and 3 (the multiple summation rule)  should be applied to the sum $\sum_{j \in \YY^+_{rr}} p_{rj}(\e)$. Second,  propositions {\bf (i)} (the multiplication by constant $-1$) and  {\bf (ii)} (the summation with constant $1$)  of Lemmas 2$_*$ and 2  should be applied to the asymptotic expansion for probability ${p}_{rr}(\e)$ given in condition  ${\bf D'}$, in order to get the asymptotic expansion for  function $1 - p_{rr}(\e)$. Third, Lemmas 1$_*$ and 1  should be applied to the asymptotic expansion for  function $\bar{p}_{rr}(\e)$ given in two alternative forms by 
relation (\ref{hopte}). This yields the corresponding pivotal the $(\bar{l}^-_{rr}, \bar{l}^+_{rr})$-expansion 
 for probabilities $\bar{p}_{rr}(\e)$,  given in Lemma 8$_*$, and proves that this expansion is a
$(\bar{l}^-_{rr}, \bar{l}^+_{rr},  \bar{\delta}_{rr}, \bar{G}_{rr}, \bar{\e}_{rr})$-expansion,  with parameters computed in the process of realization of the above algorithm.  The case $r \notin \YY_r$ is trivial, since,  in this case,  probability $\bar{p}_{rr}(\e) \equiv 1$. $\Box$
 
Let us recall formula (19)$_*$ for the transition probabilities $_rp_{ij}(\e), i, j \in \, _r\XX = \XX \setminus \{r \}$ of the reduced embedded Markov chain $_r\eta^{(\e)}_n$, introduced in Part I,  
\begin{equation}\label{transitoi}
_rp_{ij}(\e)  = p_{ij}(\e) + p_{ir}(\e) \frac{p_{rj}(\e)}{1 - p_{rr}(\e)}.
\end{equation}

Let us introduce parameter,
\begin{equation}\label{hoht}
\delta^\circ =  \min_{j \in \YY_i, i \in \XX} \delta_{ij}.
\end{equation} 

Obviously,  inequalities $\delta_{ij} \geq \delta^\circ, j \in \YY_i, i \in \XX$ hold for parameters $\delta_{ij}$ penetrating upper bounds for the remainders of asymptotic expansions in condition  ${\bf D'}$. \vspace{1mm}

{\bf Theorem 1}. {\em  Conditions ${\bf A}$ and  ${\bf D'}$, assumed to hold for the Markov chains $\eta^{(\e)}_n$, also hold for the reduced Markov chains 
$_r\eta^{(\e)}_n$, for every $r \in \XX$. Also, for every  $j \in \, _r\YY_{i}, \ i \in \, _r\XX,\,  r \in \XX$, the pivotal $(_rl_{ij}^-, \, _rl_{ij}^+)$-expansion  for the transition probability $_rp_{ij}(\e)$ given in Theorem 2$_*$ is a $(_rl_{ij}^-, \, _rl_{ij}^+, \, _r\delta_{ij}, \, _rG_{ij}$, $_r\e_{ij})$-expansion penetrating condition ${\bf D'}$ for the Markov chains $_r\eta^{(\e)}_n$. Parameters $ _r\delta_{ij}, \, _rG_{ij}$ and $_r\e_{ij}$ can be computed using the algorithm  described below,  in the proof of the theorem. The inequalities $_r\delta_{ij} \geq \delta^\circ, \, j \in \, _r\YY_{i}, \, i \in \, _r\XX, \, r \in \XX$ hold. 
} \vspace{1mm}

{\bf Proof}. Condition ${\bf A}$ holds for the Markov chains $_r\eta^{(\e)}_n$ by Lemma  6$_*$, with the same parameter $\e_0$ as for the Markov chains $\eta^{(\e)}_n$  and with the transition sets $_r\YY_{i}, i \in \, _r\XX$ given  by relation (20)$_*$.  

Let us prove that condition ${\bf D'}$ holds for the Markov chains $_r\eta^{(\e)}_n$, with the same parameter $\e_0$ as for the Markov chains $\eta^{(\e)}_n$   and the transition sets $_r\YY_{i}, i \in \, _r\XX$ given  by relation (20)$_*$. Let $j, r \in \YY_i \cap \YY_r$. First, propositions {\bf (v)} (the division rule) of Lemmas 2$_*$ and 2 should be applied to the quotient  $\frac{p_{rj}(\e)}{1 -  p_{rr}(\e)}$. Second,  propositions {\bf (iii)} (the multiplication rule) of Lemmas 2$_*$ and 2  should be applied to the product $p_{ir}(\e) \cdot 
\frac{p_{rj}(\e)}{1 -  p_{rr}(\e)}$. Third,  propositions {\bf (ii)} (the summation rule) of Lemmas 2$_*$ and 2  should be applied to sum $_rp_{ij}(\e)  = p_{ij}(\e) +  p_{ir}(\e) 
\cdot \frac{p_{rj}(\e)}{1 -  p_{rr}(\e)}$. The asymptotic expansions for probabilities $p_{ir}(\e), p_{rj}(\e)$, and $p_{ij}(\e)$, given in condition ${\bf D'}$, and probability 
$1 -  p_{rr}(\e)$, given in Lemmas 8$_*$ and 5, should be used.  This yields the corresponding pivotal 
$(_rl_{ij}^-, \, _rl_{ij}^+)$-expansions for transition probabilities  $_rp_{ij}(\e),  j \in \, _r\YY_{i}, \ i \in \, _r\XX, \, r \in \XX$,  given  in Theorem 2$_*$, and proves that these expansions are $(_rl_{ij}^-, \, _rl_{ij}^+, \, _r\delta_{ij}, \, _rG_{ij}$, $_r\e_{ij})$-expansions, with parameters computed in the process of realization  of the above algorithm. If $j \notin \YY_i$ then $p_{ij}(\e) \equiv 0$; if $j \notin \YY_r$ then $p_{rj}(\e) \equiv 0$; if $r \notin \YY_i$ then $p_{ir}(\e) \equiv 0$; if $r \notin \YY_r$ then $1 - p_{rr}(\e) \equiv 1$. In these cases,  the above algorithm is readily simplified. Thus, condition ${\bf D'}$ holds for the reduced Markov chains $_r\eta^{(\e)}_n$. 

Inequalities $_r\delta_{ij} \geq \delta^\circ, \, j \in \, _r\YY_{i}, \, i \in \, _r\XX, \, r \in \XX$ hold,  by proposition {\bf (iii)} of Lemma 4. $\Box$   

Let us recall formula (22)$_*$ for the expectations of sojourn times   $_re_{ij}(\e)$, for $i, j \in \, _r\XX = \XX \setminus \{r \}$ for the  reduced semi-Markov process  
$_r\eta^{(\e)}(t)$,  introduced in Part I,  
\begin{align}\label{expaga}
_re_{ij}(\e)   & = e_{ij}(\e)  + e_{ir}(\e) \frac{p_{rj}(\e)}{1 - p_{rr}(\e)}  \nonumber \\ 
& \quad +  e_{rr}(\e)  \frac{p_{ir}(\e)}{1 - p_{rr}(\e)} \frac{p_{rj}(\e) }{1 - p_{rr}(\e)} +  e_{rj}(\e)  \frac{p_{ir}(\e)}{1 - p_{rr}(\e)}. 
\end{align}

Let us introduce parameter,
\begin{equation}\label{parawa}
\delta^* =  \min_{j \in \YY_i, i \in \XX} ( \delta_{ij} \wedge \dot{\delta}_{ij} ).
\end{equation}

Obviously,  inequalities $\delta_{ij}, \dot{\delta}_{ij} \geq \delta^*, j \in \YY_i, i \in \XX$ hold for parameters $\delta_{ij}$ and $\dot{\delta}_{ij}$ penetrating upper bounds for the remainders of asymptotic expansions in conditions  ${\bf D'}$ and ${\bf E'}$. \vspace{1mm}

{\bf Theorem 2}. {\em  Conditions ${\bf A}$ -- ${\bf C}$, ${\bf D'}$ and ${\bf E'}$, assumed to hold for the semi-Markov processes $\eta^{(\e)}(t)$, also hold for the reduced semi-Markov processes $_r\eta^{(\e)}(t)$, for every $r \in \XX$. Also, for every $j \in \, _r\YY_{i}, \ i \in \, _r\XX, \, r \in \XX$, the pivotal 
 $(_rm_{ij}^-, \, _rm_{ij}^+)$-expansion for expectation  $_re_{ij}(\e)$  given in Theorem 3$_*$ is a  
 $(_rm_{ij}^-, \, _rm_{ij}^+, \, _r\dot{\delta}_{ij}, \, _r\dot{G}_{ij}, \, _r\dot{\e}_{ij})$-expansion penetrating condition ${\bf E'}$ for the semi-Markov processes $_r\eta^{(\e)}(t)$. Parameters  $_r\dot{\delta}_{ij}, \, _r\dot{G}_{ij}$ and $_r\dot{\e}_{ij}$ can be computed using the algorithm  described below,  in the proof of the theorem. The inequalities $_r\dot{\delta}_{ij} \geq \delta^*, \, j  \in \, _r\YY_{i}, \,  i \in \, _r\XX, \, r \in \XX$ hold. 
} \vspace{1mm}

{\bf Proof}. Conditions ${\bf A}$ and ${\bf D'}$ hold for the semi-Markov processes $_r\eta^{(\e)}(t)$, respectively, by Lemma  6$_*$ and Theorem 1, with the same parameter $\e_0$ as for the semi-Markov processes  $\eta^{(\e)}(t)$, and the transition sets $_r\YY_{i}, i \in \, _r\XX$ given  by relation (20)$_*$.  Also conditions ${\bf B}$ and ${\bf C}$ hold for processes $_r\eta^{(\e)}(t)$,  by Lemma 7$_*$.  

Let us prove that condition ${\bf E'}$ holds for the semi-Markov processes $_r\eta^{(\e)}(t)$, with the same parameter $\e_0$ and the transition sets $_r\YY_{i}, i \in \, _r\XX$ given  by relation (20)$_*$.  Let $j, r \in \YY_i \cap \YY_r$. First,  propositions {\bf (v)} (the division rule) of Lemmas 2$_*$ and 2 should be applied to the quotients  
$\frac{p_{rj}(\e)}{1 -  p_{rr}(\e)}$ and $\frac{p_{ir}(\e)}{1 -  p_{rr}(\e)}$. Second,  propositions {\bf (iii)} (the multiplication rule) of Lemmas 2$_*$ and 2  should be applied to  the products $e_{ir}(\e) \cdot \frac{p_{rj}(\e)}{1 -  p_{rr}(\e)}$ and  $e_{rj}(\e) \cdot \frac{p_{ir}(\e)}{1 -  p_{rr}(\e)}$ and propositions {\bf (ii)} (the multiple multiplication rule) of Lemmas 3$_*$ and 3  to the  product $e_{rr}(\e) \cdot  \frac{p_{ir}(\e)}{1 -  p_{rr}(\e)} \cdot \frac{p_{rj}(\e)}{1 -  p_{rr}(\e)}$.
Third,  propositions {\bf (i)} (the multiple summation rule) of Lemmas 3$_*$ and 3  should be applied to sum $_re_{ij}(\e)  = e_{ij}(\e) +  e_{ir}(\e) \cdot \frac{p_{rj}(\e)}{1 -  p_{rr}(\e)} +  e_{rr}(\e) \cdot \frac{p_{ir}(\e)}{1 -  p_{rr}(\e)} \cdot \frac{p_{rj}(\e)}{1 -  p_{rr}(\e)} + e_{rj}(\e) \cdot \frac{p_{ir}(\e)}{1 -  p_{rr}(\e)}$. The asymptotic expansions for probabilities $p_{ir}(\e), p_{ir}(\e)$ and $p_{ij}(\e)$, given in condition ${\bf D'}$, probability $1 -  p_{rr}(\e)$, given in Lemmas 8$_*$ and 5, and expectations $e_{ij}(\e), e_{ir}(\e), e_{rr}(\e)$ and  $e_{rj}(\e)$, given in condition ${\bf E'}$, should be used.  This, first, yields the corresponding pivotal $(_rm_{ij}^-, \, _rm_{ij}^+)$-expansions 
for expectations of sojourn times $_re_{ij}(\e),  j \in \, _r\YY_{i}, \ i \in \, _r\XX, \, r \in \XX$, s given  in Theorem 3$_*$, and, second, proves that these expansions are $(_rm_{ij}^-, \, _rm_{ij}^+, \, _r\dot{\delta}_{ij}, \, _r\dot{G}_{ij}, \, _r\dot{\e}_{ij})$-expansions,  with parameters computed in the process of realization  of the above algorithm. 
If $j \notin \YY_i$ then $p_{ij}(\e) \equiv 0$ and $e_{ij}(\e) \equiv 0$; if $j \notin \YY_r$ then $p_{rj}(\e) \equiv 0$ and $e_{rj}(\e) \equiv 0$; if $r \notin \YY_i$ then $p_{ir}(\e) \equiv 0$ and $e_{ir}(\e) \equiv 0$; if $r \notin \YY_r$ then $1 - p_{rr}(\e) \equiv 1$ and $e_{rr}(\e) \equiv 0$. In these cases, the above algorithm is readily simplified. 
Thus, condition ${\bf E'}$ holds for the reduced semi-Markov processes  $_r\eta^{(\e)}(t)$. 

Inequalities $_r\dot{\delta}_{ij} \geq \delta^*, \, j \in \, _r\YY_{i}, \, i \in \, _r\XX, \, r \in \XX$ hold,  by proposition {\bf (iii)} of Lemma 4. $\Box$

It worth to note that,  despite bulky forms, formulas for parameters of upper bounds for remainders, in  the asymptotic expansions given in Lemma 5 and Theorems 1 and 2, are computationally effective. \\

{\bf 4. Asymptotic expansions for expected hitting times \\ \makebox[11mm]{} with explicit upper bounds for remainders}  \\

As in Part I,  let $\bar{r}_{i, N}  = \langle r_{i, 1}, \ldots, r_{i, N} \rangle =  \langle r_{i, 1}, \ldots, r_{i, N-1}, i \rangle$ be  a permutation of the sequence  $\langle 1, \ldots, N \rangle$ such that 
$r_{i, N} = i$, and let $\bar{r}_{i, n} = \langle r_{i, 1}, \ldots, r_{i, n} \rangle$, $n = 1, \ldots, N$ be the corresponding chain of growing sequences of states from space $\XX$. \vspace{1mm}

{\bf Theorem 3}. {\em Let conditions ${\bf A}$ -- ${\bf C}$, ${\bf D'}$ and  ${\bf E'}$ hold for the semi-Markov processes  $\eta^{(\e)}(t)$. Then,  for every $i \in \XX$, the pivotal $(M_{ii}^-, M_{ii}^+)$-expansion for  the expectation of hitting time $E_{ii}(\e)$, given in Theorem 4$_*$ and obtained as the result of sequential exclusion of states  $r_{i, 1}, \ldots, r_{i, N-1}$ from  the phase space $\XX$  of the processes  $\eta^{(\e)}(t)$,  is a $(M_{ii}^-, M_{ii}^+,  \, _{\bar{r}_{i, N - 1}}\dot{\delta}_{ii}, \, _{\bar{r}_{i, N - 1}}\dot{G}_{ii}, \, _{\bar{r}_{i, N - 1}}\dot{\e}_{ii})$-expansion. Parameters $_{\bar{r}_{i, N - 1}}\dot{\delta}_{ii},  \, _{\bar{r}_{i, N - 1}}\dot{G}_{ii}$ and  $_{\bar{r}_{i, N - 1}}\dot{\e}_{ii}$ can be computed using the algorithm  described below,  in the proof of the theorem. Also, inequality $_{\bar{r}_{i, N - 1}}\dot{\delta}_{ii} \geq \delta^*$ holds making it possible to rewrite   function $E_{ii}(\e)$ as the pivotal $(M_{ii}^{-}, M_{ii}^{+}, \delta^*, \, _{\bar{r}_{i, N - 1}}G^*_{ii}, \, _{\bar{r}_{i, N - 1}}\dot{\e}_{ii})$-expansion, with parameter $_{\bar{r}_{i, N - 1}}G^*_{ii} = \, _{\bar{r}_{i, N - 1}}\dot{G}_{ii} \,  
\cdot (_{\bar{r}_{i, N - 1}} \dot{\e}_{ii})^{(_{\bar{r}_{i, N - 1}}\dot{\delta}_{ii} - \delta^*)}$. } \vspace{1mm}

{\bf Proof}.  Let us assume that $p^{(\e)}_i = 1$.  Denote as $_{\bar{r}_{i, 0}}\eta^{(\e)}(t) = \eta^{(\e)}(t)$ the initial semi-Markov process. Let us exclude  state 
$r_{i, 1}$ from the phase space of semi-Markov process $_{\bar{r}_{i, 0}}\eta^{(\e)}(t)$ using the time-space screening procedure described in Section 5$_*$. Let $_{\bar{r}_{i, 1}}\eta^{(\e)}(t)$ be the corresponding reduced semi-Markov process. The above procedure can be repeated. The state $r_{i, 2}$ can  be excluded from the phase space of the semi-Markov process $_{\bar{r}_{i, 1}}\eta^{(\e)}(t)$. Let $_{\bar{r}_{i, 2}}\eta^{(\e)}(t)$ be the corresponding reduced semi-Markov process. By continuing the above procedure for states $r_{i, 3}, \ldots, r_{i, n}$,  we construct the reduced semi-Markov process $_{\bar{r}_{i, n}}\eta^{(\e)}(t)$. 

The process   $_{\bar{r}_{i, n}}\eta^{(\e)}(t)$ has  the phase space $_{\bar{r}_{i, n}}\XX = \XX \setminus \{ r_{i, 1}, r_{i, 2}, \ldots, r_{i, n} \}$. 
The transition probabilities  $_{\bar{r}_{i, n}}p_{i'j'}(\e)$,  
$i', j' \in \, _{\bar{r}_{i, n}}\XX$ and the expectations of sojourn times $_{\bar{r}_{i, n}}e_{i' j'}(\e), i', j' \in \, _{\bar{r}_{i, n}}\XX$  are determined for the process $_{\bar{r}_{i, n}}\eta^{(\e)}(t)$  by the transition probabilities  and the expectations of sojourn times for the process $_{\bar{r}_{i, n-1}}\eta^{(\e)}(t)$,  via relations (\ref{transitoi}) and (\ref{expaga}). 

By Theorem 1$_*$, the expectation of hitting time $E_{i'j'}(\e)$ coincides for the semi-Markov processes $_{\bar{r}_{i, 0}}\eta^{(\e)}(t)$,  $_{\bar{r}_{i, 1}}\eta^{(\e)}(t), \ldots$, $_{\bar{r}_{i, n}}\eta^{(\e)}(t)$, for every  $i', j' \in \, _{\bar{r}_{i, n}}\XX$.  

By Theorems 2$_*$, 3$_*$,  1 and 2, the  semi-Markov processes $_{\bar{r}_{i, n}}\eta^{(\e)}(t)$ satisfy conditions ${\bf A}$ -- ${\bf C}$, ${\bf D'}$ and  ${\bf E'}$. The transition sets $_{\bar{r}_{i, n}}\YY_{i'}, i' \in \, _{\bar{r}_{i, n}}\XX$, for the process $_{\bar{r}_{i, n}}\eta^{(\e)}(t)$, are determined  by  the transition sets $_{\bar{r}_{i, n-1}}\YY_{i'}, i' \in \, _{\bar{r}_{i, n-1}}\XX$, for the process $_{\bar{r}_{i, n-1}}\eta^{(\e)}(t)$, via 
relation (20)$_*$.  For every $j' \in \,  _{\bar{r}_{i, n}}\YY_{i'}, i' \in \,  _{\bar{r}_{i, n}}\XX$, the pivotal $(_{\bar{r}_{i, n}}l_{i'j'}^-, \, _{\bar{r}_{i, n}}l^+_{i'j'})$-expansion for  transition probability $_{\bar{r}_{i, n}}p_{i'j'}(\e)$, given in Theorem 2$_*$, is, by Theorem 1, a $(_{\bar{r}_{i, n}}l_{i'j'}^-, \, _{\bar{r}_{i, n}}l^+_{i'j'}, \, _{\bar{r}_{i, n}}\delta_{i'j'}, \, _{\bar{r}_{i, n}}G_{i'j'}$, $_{\bar{r}_{i, n}}\e_{i'j'})$-expansion, with parameters $_{\bar{r}_{i, n}}\delta_{i'j'}$,  $_{\bar{r}_{i, n}}G_{i'j'}$ and $_{\bar{r}_{i, n}}\e_{i'j'}$ given in this theorem. Analogously, for every $j' \in \,  _{\bar{r}_{i, n}}\YY_{i'}, i' \in \,  _{\bar{r}_{i, n}}\XX$, the pivotal 
$(_{\bar{r}_{i, n}}m_{i'j'}^-$, $_{\bar{r}_{i, n}}m_{i'j'}^+)$-expansion for expectation $_{\bar{r}_{i, n}}e_{i'j'}(\e)$, given in Theorem 3$_*$, is, by Theorem 2, a $(_{\bar{r}_{i, n}}m_{i'j'}^-, \, _{\bar{r}_{i, n}}m_{i'j'}^+,  \, _{\bar{r}_{i, n}}\dot{\delta}_{i'j'}$,  
$_{\bar{r}_{i, n}}\dot{G}_{i'j'}, _{\bar{r}_{i, n}}\dot{\e}_{i'j'})$-expansion, with parameters  $_{\bar{r}_{i, n}}\dot{\delta}_{i'j'}, \, _{\bar{r}_{i, n}}\dot{G}_{i'j'}$ and  $_{\bar{r}_{i, n}}\dot{\e}_{i'j'}$ given in this theorem. Also, by Theorem 2, the inequalities 
$_{\bar{r}_{i, n}}\dot{\delta}_{i' j'} \geq \delta^*, j' \in \,  _{\bar{r}_{i, n}}\YY_{i'}, i' \in \,  _{\bar{r}_{i, n}}\XX$ hold.

The algorithm described above has a recurrent form and should be realized sequentially for the reduced semi-Markov processes \, $_{\bar{r}_{i, 1}}\eta^{(\e)}(t), \ldots$, 
$_{\bar{r}_{i, n}}\eta^{(\e)}(t)$ starting from the initial semi-Markov process $ _{\bar{r}_{i, 0}}\eta^{(\e)}(t)$.

Let us take $n = N-1$. The semi-Markov process  $_{\bar{r}_{i, N-1}}\eta^{(\e)}(t)$ has the phase space $_{\bar{r}_{i, N-1}}\XX = \XX \setminus \{ r_{i, 1}, r_{i, 2}, \ldots, r_{i, N- 1} \} = \{ i \}$,  which is a one-state set. The process  
$_{\bar{r}_{i, N-1}}\eta^{(\e)}(t)$ returns to state $i$ after every jump. Its transition probability $_{\bar{r}_{i, N-1}}p_{ii}(\e) = 1$,  and the expectation of hitting time,  $E_{ii}(\e) = \, _{\bar{r}_{i, N-1}}e_{ii}(\e)$. This equality and Theorem 4$_*$ yield, for every $i \in \XX$, the  pivotal $(M_{ii}^{-}, M_{ii}^{+})$-expansion for expectation $E_{ii}(\e)$, which is invariant with respect to  any permutation $\bar{r}'_{i, N-1} =  \langle r'_{i, 1}, \ldots, r'_{i, N-1} \rangle $ of the sequence $\bar{r}_{i, N-1} =  \langle r_{i, 1}, \ldots, r_{i, N-1} \rangle$. This invariance also implies that parameters $_{\bar{r}_{i, N-1}}m_{ii}^\pm 
= M_{ii}^{\pm}$ do not depend on the choice of sequence $\bar{r}_{i, N-1}$. The above $(M_{ii}^{-}, M_{ii}^{+})$-expansion for  $E_{ii}(\e) \ = \ _{\bar{r}_{i, N-1}}e_{ii}(\e)$ is  a $(M_{ii}^-, M_{ii}^+,  \ _{\bar{r}_{i, N-1}}\dot{\delta}_{ii}, \ _{\bar{r}_{i, N-1}}\dot{G}_{ii}$, $_{\bar{r}_{i, N-1}}\dot{\e}_{ii})$-expansion.

The above algorithm can be realized for any sequence $\bar{r}_{i, N -1}  = \langle r_{i, 1}, \ldots$, $r_{i, N-1} \rangle$, but the invariance  of explicit 
upper bounds for  remainders, with respect to permutations  $\bar{r}_{i, N}  = \langle r_{i, 1}, \ldots, r_{i, N-1}, i \rangle$ of sequence $\langle 1, \ldots, N \rangle$,  can not be guaranteed. 

However, the inequality 
$_{\bar{r}_{i, N-1}}\dot{\delta}_{ii} \geq \delta^*$ holds, for any sequence $\bar{r}_{i, N-1} =  \langle r_{i, 1}, \ldots,$ $r_{i, N-1} \rangle$, by Theorem 2.

The algorithm described above can  be repeated, for every $i \in \XX$. $\Box$

It is worth to note that the algorithms based on sequential exclusion of states from the phase space of perturbed semi-Markov processes make it possible to get Laurent asymptotic expansions (without and with explicit upper bounds for remainders) for expected hitting times, for nonlinearly perturbed semi-Markov processes. Such asymptotic results have their own important value.

Let $\bar{r}_{i, j, N}  = \langle r_{i, j, 1}, \ldots, r_{i, j, N} \rangle =  \langle r_{i, j, 1}, \ldots, r_{i, j, N-2}, i , j\rangle$ be  a permutation of the sequence  $\langle 1, \ldots, N \rangle$ such that 
$r_{i, j, N -1} = i, r_{i, j, N} = j$, and let  $\bar{r}_{i, j, n} = \langle r_{i, j, 1}, \ldots, r_{i, j, n} \rangle$, $n = 1, \ldots, N$ be the corresponding chain of growing sequences of states from space $\XX$. 

By applying the algorithm of sequential phase space reduction described in Theorem 4 to the above sequence of states  
$\bar{r}_{i, j, N-2}$, we construct the reduced semi-Markov process $_{\bar{r}_{i, j, N-2}}\eta^{(\e)}(t)$. This process has the phase space 
$_{\bar{r}_{i, j, N-2}}\XX = \XX_{ij} = \{i, j \}$, which is a two-states set. The transition probabilities of the embedded Markov chain 
$_{\bar{r}_{i, j, N-2}}p_{i'j'}(\e) = p_{ij, i'j'}(\e), i', j' \in \XX_{ij}$, the expectations of sojourn times $_{\bar{r}_{i, j, N-2}}e_{i'j'}(\e)$ $= e_{ij, i'j'}(\e), i', j'$ $\in \XX_{ij}$,  and the transition sets $_{\bar{r}_{i, j, N-2}}\YY_{i'} = \YY_{ij, i'}, i' \in \XX_{ij}$ can be found using the recurrent algorithm described in Theorem 4$_*$. These probabilities, expectations and transition sets are invariant to any permutation  $\bar{r}'_{i, j, N-2}$ of sequence $\bar{r}_{i, j, N-2}$. This legitimates the above alternative simplified notations.  

Theorem 4$_*$ let us construct the pivotal  $(l_{ij, i'j'}^-, l_{ij, i'j'}^+)$-expansions for the transition probabilities   $p_{ij, i'j'}(\e), j' \in \YY_{ij, i'},   i' \in \XX_{ij}$  and the pivotal $(m_{ij, i'j'}^-$, $m_{ij, i'j'}^+)$-expansions for the expectations of hitting times   $e_{ij, i'j'}(\e),  j' \in \YY_{ij, i'}$, $i' \in \XX_{ij}$,  using the recurrent algorithm based on sequential exclusion states 
$r_{i, j, 1}, \ldots, r_{i, j, N-2}$ from the phase space $\XX$. These expansions are invariant to any permutation  $\bar{r}'_{i, j, N-2}$ of sequence $\bar{r}_{i, j, N-2}$.  

According Theorem 3,   the above  $(l_{ij, i'j'}^-, l_{ij, i'j'}^+)$-expansions for transition probabilities  are   $(l_{ij, i'j'}^-, l_{ij, i'j'}^+, \ _{\bar{r}_{i, j, N-2}}\delta_{ i'j'}, \ _{\bar{r}_{i, j, N-2}}G_{ i'j'}, \, _{\bar{r}_{i, j, N-2}}\e_{ i'j'})$-expansions and  the above  $(m_{ij, i'j'}^-, m_{ij, i'j'}^+)$-expansions for  expectations of sojourn times   are $(m_{ij, i'j'}^-, \, m_{ij, i'j'}^+$, $_{\bar{r}_{i, j, N-2}}\dot{\delta}_{ i'j'}, _{\bar{r}_{i, j, N-2}}\dot{G}_{ i'j'}$, $_{\bar{r}_{i, j, N-2}}\dot{\e}_{ i'j'})$-expansions. This theorem let us also  compute  parameters  $_{\bar{r}_{i, j, N-2}}\delta_{ i'j'}, \, _{\bar{r}_{i, j, N-2}}G_{ i'j'}, \, _{\bar{r}_{i, j, N-2}}\e_{ i'j'}$  and $_{\bar{r}_{i, j, N-2}}\dot{\delta}_{ i'j'}, \, _{\bar{r}_{i, j, N-2}}\dot{G}_{ i'j'}, \, _{\bar{r}_{i, j, N-2}}\dot{\e}_{ i'j'}$ of upper bounds for the corresponding remainders.  

By Theorem 1$_*$, the expectation of hitting time $E_{i', j'}(\e)$ coincides for the initial semi-Markov processes $\eta^{(\e)}(t)$ and the reduced semi-Markov process $_{\bar{r}_{i, j, N-2}}\eta^{(\e)}(t)$, for every $i', j' \in \XX_{ij}$.  This obviously implies that these expectations are also invariant to any permutation  $\bar{r}'_{i, j, N-2}$ of sequence $\bar{r}_{i, j, N-2}$.

It is easy to write down the formulas for the above  expectations,  for the two-states semi-Markov process $_{\bar{r}_{i, j, N-2}}\eta^{(\e)}(t)$. These formulas are,  
$E_{i'j'}(\e)  = e_{ij, i'}(\e)  \frac{1}{1 - p_{ij, i'i'}(\e)}$,  $E_{j'j'}(\e)  = e_{ij, j'}(\e) \, + \,  e_{ij, i'}(\e) \frac{p_{ij, j'i'}(\e) }{1 - p_{ij, i'i'}(\e)}$, 
where $e_{ij, i'}(\e) =$  $e_{ij, i'i'}(\e) + e_{ij, i' j'}(\e), \ i', j' \in \XX_{ij}, i' \neq j'$.

Under the assumption that conditions of Theorem 4$_*$ hold, the operational rules given in Lemma 2$_*$ can be applied to functions $E_{i'j'}(\e), i', j' \in \XX_{ij}$, in order to get the corresponding $(M^-_{i'j'}, M^+_{i'j'})$-expansions.  These expansions are  invariant to any permutation  $\bar{r}'_{i, j, N-2}$ of sequence $\bar{r}_{i, j, N-2}$, used 
in the corresponding recurrent algorithm based on sequential exclusion states $r_{i, j, 1}, \ldots, r_{i, j, N-2}$ from the phase space $\XX$.

Finally, under the assumption that conditions of Theorem 3 hold,  the operational rules given in Lemma 2 can  be applied, in order to prove that the above \, $(M^-_{i'j'}, \, M^+_{i'j'})$-expansions are $(M^-_{i'j'}, \, M^+_{i'j'}, \ {_{\bar{r}_{i, j N-2}}}\dot{\delta}_{i'j'}, \ {_{\bar{r}_{i, j, N-2}}}\dot{G}_{i'j'}$, ${_{\bar{r}_{i, j, N-2}}}\dot{\e}_{i'j'}$)-expansions,  and to compute parameters  $_{\bar{r}_{i, j N-2}}\dot{\delta}_{i'j'}$, $_{\bar{r}_{i, j, N-2}}\dot{G}_{i'j'}$ and  $_{\bar{r}_{i, j, N-2}}\dot{\e}_{i'j'}$ of upper bounds for the corresponding remainders. Also, by Lemma 4, the inequality  $_{\bar{r}_{i, j N-2}}\dot{\delta}_{i'j'} \geq \delta^*$ holds, for every $i', j' \in \XX_{ij}$, sequence $\bar{r}_{i, j N-2}$, and $i, j \in \XX$. \\

{\bf 5. Asymptotic expansions for stationary distributions with \\ \makebox[10mm]{} explicit upper bounds for remainders}  \\

Let us recall the pivotal $(n_{i}^-, n_{i}^+)$-expansion for stationary probability $\pi_{i}(\e) $ of nonlinearly perturbed semi-Markov process 
$\eta^{(\e)}(t)$ given, under conditions ${\bf A}$ -- ${\bf E}$, in Theorem 5$_*$. This asymptotic expansion has the following form,   for $i \in \XX$, 
\begin{equation}\label{expanaba}
\pi_{i}(\e) =  \sum_{l = n_{i}^-}^{n_{i}^+} c_{i}[l]\e^l + o_i(\e^{n_{i}^+}), \, \e \in (0, \e_0].
\end{equation}
 
According Theorem 5$_*$, the above asymptotic expansion  is invariant with respect to the choice of sequence states $\bar{r}_{i, N-1} = (r_{i, 1}, \ldots, r_{i, N-1})$ used in the corresponding algorithm, for every $i \in \XX$.  

The following theorem is the main new result in Part II of this paper. \vspace{1mm}

{\bf Theorem 4}. {\em  Let conditions ${\bf A}$ -- ${\bf C}$, ${\bf D'}$ and  ${\bf E'}$ hold for the semi-Markov processes  $\eta^{(\e)}(t)$.  Then,  for every $i \in \XX$, the pivotal $(n_{i}^-, n_{i}^+)$-expansion {\rm (\ref{expanaba})} for the stationary probability $\pi_{i}(\e)$, given  in Theorem 5$_*$ and obtained as the result of sequential exclusion of states $r_{i, 1}, \ldots, r_{i, N-1}$ from  the phase space $\XX$  of the processes  $\eta^{(\e)}(t)$, is a 
$(n_{i}^-, n_{i}^+, \, _{\bar{r}_{i, N - 1}}\delta_i, \, _{\bar{r}_{i, N - 1}}G_i, \, _{\bar{r}_{i, N - 1}}\e_i)$-expansion.  Parameters \,  $_{\bar{r}_{i, N - 1}}\delta_i, \, _{\bar{r}_{i, N - 1}}G_i$ and $_{\bar{r}_{i, N - 1}}\e_i$ can be computed using the algorithm  described below,  in the proof of the theorem.
Also, inequality $_{\bar{r}_{i, N - 1}}\delta^*_{i} \geq \delta^*$ holds making it possible to rewrite   function $\pi_{i}(\e)$ as the pivotal $(n_{i}^{-}, n_{i}^{+}, \delta^*, \, _{\bar{r}_{i, N - 1}}G^{*}_{i}, \, _{\bar{r}_{i, N - 1}}\e_{i})$-expansion, with parameter $_{\bar{r}_{i, N - 1}}G^{*}_{i} = \, _{\bar{r}_{i, N - 1}}G_{i} \cdot 
\, (_{\bar{r}_{i, N - 1}} \e_{i})^{(_{\bar{r}_{i, N - 1}}\delta_{i} - \delta^*)}$. } \vspace{1mm}

{\bf Proof}. Let us choose an arbitrary state $i \in \XX$. First, proposition {\bf (i)} (the multiple summation rule) of Lemmas 3$_*$ and 3 should be applied to the pivotal $(m_{i}^-, m_{i}^+)$-expansion for the expectation $e_{i}(\e)  = \sum_{j \in \YY_i} e_{ij}(\e)$ given by  relation (29)$_*$, in the proof of Theorem 5$_*$. This yields a $(m_{i}^-, m_{i}^+, 
\dot{\delta}_i, \dot{G}_i, \dot{\e}_i)$-expansion for the expectation $e_{i}(\e)$, with the corresponding parameters 
$\dot{\delta}_i, \dot{G}_i$ and  $\dot{\e}^*_i$. Second,  
the propositions {\bf (v)} (the division rule) of Lemmas 2$_*$ and 2 should be applied to the quotient $\pi_i(\e) = \frac{e_i(\e)}{E_{ii}(\e)}$. The $(m_{i}^-, m_{i}^+, \dot{\delta}_i, \dot{G}_i, \dot{\e}_i)$-expansion for the expectation $e_{i}(\e)$ and the $(M_{ii}^-, \, M_{ii}^+,  \ _{\bar{r}_{i, N-1}}\dot{\delta}_{ii}, \ _{\bar{r}_{i, N-1}}\dot{G}_{ii}$, $_{\bar{r}_{i, N-1}}\dot{\e}_{ii})$-expansion for the expectation of hitting time $E_{ii}(\e)$, given in Theorems 4$_*$ and 3, should be used.  This yields the corresponding pivotal 
 $(n_{i}^-, n_{i}^+)$-expansion for stationary probability $\pi_i(\e)$, given in Theorem 5$_*$, and proves that
this expansion is a $(n_{i}^-, n_{i}^+, \, _{\bar{r}_{i, N - 1}}\delta_i, \, _{\bar{r}_{i, N - 1}}G_i, \, _{\bar{r}_{i, N - 1}}\e_i)$-expansion,  with parameters computed in the process of realization  of the above algorithm. 
Inequality $_{\bar{r}_{i, N - 1}}\delta_{i} \geq \delta^*$ holds, for every sequence $\bar{r}_{i, N - 1}$, by proposition {\bf (iii)} of Lemma 4. $\Box$

The explicit upper bounds for remainders in the asymptotic expansions  given in Theorem 4 have a clear and informative power-type form. An  useful property of these upper bounds  is that they are uniform with respect to the perturbation parameter. The recurrent algorithm  for finding these upper bounds is computationally effective. 

Unfortunately,  the summation and multiplication operational rules for computing power-type upper bounds for remainders  possess commutative but do not possess associative and distributive properties. This causes 
dependence of the resulting upper bounds for remainders in the asymptotic expansions for stationary probabilities $\pi_i(\e), \, i \in \XX$ on a choice of the corresponding sequences of states
 $\bar{r}_{i, N-1} =  \langle r_{i, 1}, \ldots, r_{i, N-1} \rangle,  \, i \in \XX$ used in the above algorithm. This rises two open questions, the first one, about possible  alternative forms for remainders possessing the desirable algebraic properties mentioned above, and, the second one,  about an optimal choice of sequences of states $\bar{r}_{i, N-1},  \, i \in \XX$. 
 
In conclusion,  we  would like to mention some prospective directions for future research studies. 

The method of sequential reduction of phase space presented in the paper  can be applied  for getting asymptotic expansions for high order power  and exponential moments of hitting times,  for nonlinearly perturbed semi-Markov processes. This is an interesting problem, which has  its own important theoretical and applied values.

We are quite sure that a combination of results in the above direction with the methods of asymptotic analysis for nonlinearly perturbed regenerative processes  developed and throughly presented in  Gyllenberg and Silvestrov (2008) will make it possible to expand results from this book,  related to  asymptotic expansions for stationary and more general quasi-stationary distributions as well as other characteristics for nonlinearly perturbed semi-Markov processes with absorption,  to nonlinearly perturbed semi-Markov processes with an arbitrary asymptotic communicative structure of phase spaces.

The problems of aggregation of steps in the time-space screening procedures for semi-Markov processes, tracing  pivotal orders for different groups of states as well as  getting explicit  formulas, for coefficients and parameters of upper bounds for remainders in the corresponding asymptotic expansions for stationary distributions and moments of hitting times,  do require additional studies.  It can be expected that such formulas can be obtained, for example, for nonlinearly perturbed birth-death-type  semi-Markov processes, for which the proposed algorithms of phase space reduction   preserve the birth-death structure for reduced semi-Markov processes.

Applications  to control and queuing systems, information networks, epidemic models and  models of mathematical genetics and population dynamics, analogous to those presented in the books cited in the introduction, also create a prospective area for future research based on the asymptotic results obtained in the present paper.

\begin{center}
\bf Appendix A: Operational rules for Laurent asymptotic expansions 
\end{center}

Let us give short proofs of Lemmas 1$_*$ -- 4$_*$ and 1 -- 4 omitting some known or obvious details. 

{\bf A.1}. The formulas given in Lemmas 1$_*$ and 1 are quite obvious. 

{\bf A.2}.  The same relates to formulas  in propositions {\bf (i)} (the multiplication by a constant rule) of 
Lemmas 2$_*$  and 2. 

Proposition {\bf (ii)} (the summation rules) of Lemmas 2$_*$  and 2 can be obtained by simple accumulation of  coefficients for different powers of $\e$ and terms accumulated in the corresponding remainders, and, then,  by using obvious upper bounds for absolute values of sums of terms accumulated in the corresponding remainders.  

Proposition {\bf (iii)} (the multiplication  rule) of Lemma 2$_*$ can be proved by  multiplication of the corresponding asymptotic expansions $A(\e)$ and  $B(\e)$  and accumulation of coefficients  for 
powers  $\e^l$ for $l = h_C, \ldots, k_C$ in their product, 
\begin{equation*}
\begin{aligned}
C(\e)  & = A(\e) B(\e) \nonumber  \\
& = (a_{h_A}\e^{h_A} + \cdots + a_{k_A}\e^{k_A} + o_A(\e^{k_A})) \nonumber  \\ 
& \quad \times (b_{h_B}\e^{h_B} +
\cdots + b_{k_B}\e^{k_B} + o_B(\e^{k_B}) )  \nonumber  \\
\end{aligned}
\end{equation*}
\begin{align} \label{boply}
& = \sum_{h_C \leq l \leq k_C} \, \sum_{i + j = l, h_A \leq i \leq k_A, h_B \leq j \leq k_B } a_i b_j \e^l 
\nonumber  \\
& \quad + \sum_{k_C < i + j, h_A \leq i \leq k_A, h_B \leq j \leq k_B } a_i b_j \e^{i + j}  + \sum_{h_B \leq j \leq k_B}  b_j \e^{j} o_A(\e^{k_A}) \nonumber \\
& \quad +  \sum_{h_A \leq i \leq k_A}  a_i \e^{i} o_B(\e^{k_B}) + o_A(\e^{k_A}) o_B(\e^{k_B}) \nonumber \\
& = \sum_{h_C \leq l \leq k_C} c_{l} \e^l + o_C(\e^{k_C}),
\end{align}
where
\begin{align}\label{hopret}
 o_C(\e^{k_C}) & = \sum_{k_C < i + j, h_A \leq i \leq k_A, h_B \leq j \leq k_B } a_i b_j \e^{i + j}  + \sum_{h_B \leq j \leq k_B}  b_j \e^{j} o_A(\e^{k_A}) \nonumber \\
& \ \ +  \sum_{h_A \leq i \leq k_A}  a_i \e^{i} o_B(\e^{k_B}) + o_A(\e^{k_A}) o_B(\e^{k_B}). 
\end{align}  

Obviously,
$
\frac{o_C(\e^{k_C})}{\e^{ k_C} }  \to 0 \ {\rm as} \ \e \to 0.
$  
It should be noted that the accumulation of coefficients for powers $\e^l$ can be made in (\ref{boply}) only up to the maximal value $l = k_C = (k_A + h_B) \wedge (k_B + h_A)$,  because of  the presence in the  expression for  remainder $o_C(\e^{k_C})$ terms  $b_{h_B} \e^{h_B} o_A(\e^{k_A})$ and 
$a_{h_A} \e^{h_A} o_B(\e^{k_B})$. 

Also, relation  (\ref{hopret}) readily implies relations {\bf (a)} -- {\bf (c)}, which  determines parameters $\delta_C, G_C,  \e_C$ in proposition {\bf (iii)} of Lemma 2.

The  assumptions of proposition  {\bf (iv)} in  Lemma  2$_*$  imply that 
$\e^{- h_B} B(\e) \to  b_{h_B} \neq 0$ as $\e \to 0$. 
This relation implies that there exists $0 < \e'_0 \leq \e_0$ such that $B(\e) \neq 0$ for $\e \in (0, \e'_0]$, and, 
thus, function $C(\e) =  \frac{1}{B(\e)}$ is well defined for $\e \in (0, \e'_0]$.

Note that $h_B \leq k_B$. The  assumptions of proposition  {\bf (iv)} of Lemma 2$_*$  imply that,  $\e^{h_B} C(\e)  = (b_{h_B} + \cdots + b_{k_B}\e^{k_B - h_B} +   o_B(\e^{h_B})\e^{- h_B})^{-1}  \to b_{h_B}^{-1}  =  c_{h_C}$ as $\e \to 0$.
This relation means that function  $\e^{h_B} C(\e)$ can be represented in the form  $\e^{h_B} C(\e) = c_{h_C} + o(1)$, where $c_{h_C} = b_{h_B}^{-1}$,  or, equivalently, that 
the following representation takes place, 
$C(\e) =  c_{h_C} \e^{- h_B} + o_1(\e^{-h_B}), \e \in (0, \e'_0]$,  
where $\frac{o_1(\e^{- h_B})}{\e^{- h_B}} \to 0$  as $\e \to 0$.
 
 The latter two relations   prove proposition  {\bf (iv)} of Lemma 2$_*$, for the case $h _B = k_B$. Indeed, these relations mean that function  $C(\e) =  \frac{1}{B(\e)} $ can be represented in the form of   $(h_C, k_C)$-expansion with parameters $h_C = - h_B$, $k_C = k_B - 2h_B = - h_B = h_C$ and coefficient $c_{h_C} = b_{h_B}^{-1}$. Moreover, since $B(\e) \cdot C(\e) \equiv 1, 0 < \e \leq \e'_0$, remainder $c_1(\e)$ can be found from the following relation, $(b_{h_B} \e^{h_B} + o(\e^{h_B}))(c_{h_C} \e^{- h_B} + o_1(\e^{-h_B})) \equiv 1$ that yields  formula,  
 $o_1(\e^{-h_B}) = - \frac{c_{h_C} \e^{- h_B} o_B(\e^{h_B})}{b_{h_B}\e^{h_B} + o_B(\e^{h_B})}$. This is formula {\bf (c)} from  proposition  {\bf (iv)} of Lemma 2$_*$, for the case $h _B = k_B$. Note that, in the case $h _B = k_B$, the above asymptotic expansion  for function $C(\e)$ can not  be extended. Indeed,
$\e^{ h_B  - 1} o_1(\e^{-h_B})   = \e^{h_B -1} ( C(\e) -  c_{h_C} \e^{- h_B})  = - \frac{c_{h_C}}{b_{h_B} + o_B(\e^{h_B})\e^{-h_B}} \frac{o_B(\e^{h_B})\e^{- h_B}}{ \e}$.
The term  $\frac{o_B(\e^{h_B})\e^{- h_B}}{ \e}$ on the right hand side in the latter relation  has  an uncertain asymptotic behavior as $\e \to 0$.

Let us now assume that $h _B + 1 \leq k_B$.  
In this case, the  assumptions of proposition  {\bf (iv)} of Lemma 2$_*$  and the above asymptotic relations  imply that
$\e^{h_B  -1} o_1(\e^{-h_B})   = \e^{ h_B -1} ( C(\e) -  c_{h_C} \e^{- h_B})   = 
(b_{h_B} + \cdots + b_{k_B} \e^{k_B - h_B} + o_B(\e^{k_B}) \e^{- h_B})^{-1}  (- b_{h_B +1} c_{h_C} - \cdots$ 
$-  \, b_{k_B}c_{h_C}\e^{k_B - h_B -1}  -  o_B(\e^{k_B}) c_{h_C} \e^{- h_B -1}) 
\to  - b_{h_B}^{-1} b_{h_B +1} c_{h_C} = c_{h_C +1}$ as $\e \to 0$.  
This relation means that function  $\e^{h_B  -1}$ $\cdot o_1(\e^{-h_B})$ can be represented in the form  $\e^{h_B  -1} o_1(\e^{-h_B}) = c_{h_C +1} + o(1)$, where $c_{h_C +1} = -  b_{h_B}^{-1} b_{h_B +1} c_{h_C}$,   or, equivalently, that the following representation takes place, $C(\e) =  c_{h_C} \e^{- h_B} + c_{h_C +1} \e^{- h_B +1} + o_2(\e^{-h_B +1})$, $\e \in (0, \e'_0]$, where $\frac{o_2(\e^{-h_B +1})}{ \e^{- h_B +1}} \to 0$  as $\e \to 0$. 

The latter two relations prove proposition  {\bf (iv)} of Lemmas 2$_*$, for the case  
$h _B + 1 = k_B$. Indeed, these relations mean that function $C(\e)$ can be represented in the form of   $(h_C, k_C)$-expansion with parameters $h_C = - h_B$, $k_C = k_B - 2h_B = - h_B + 1 = h_C +1$ and coefficients $c_{h_C} = b_{h_B}^{-1}, \, c_{h_C +1} =
-  b_{h_B}^{-1} b_{h_B +1} c_{h_C}$. Moreover, since $B(\e) \cdot C(\e)  \equiv 1$, the remainder $o_2(\e^{-h_B +1})$ can be found from the following relation, $(b_{h_B} \e^{h_B} + b_{h_B +1} \e^{h_B +1} + 
 o(\e^{h_B +1}))$ $\cdot (c_{h_C} \e^{- h_B} + c_{h_C +1} \e^{- h_B +1}  + o_2(\e^{-h_B +1})) \equiv 1$ that yields  formula, $o_2(\e^{-h_B +1})
=$  $- \frac{b_{h_B +1}c_{h_C +1}\e^2  +  (c_{h_C} \e^{- h_B} + c_{h_C + 1} \e^{- h_B +1}) o_B(\e^{h_B +1})}{b_{h_B}\e^{h_B} + o_B(\e^{h_B})}$. This is formula {\bf (c)} from  proposition  {\bf (iv)} of Lemma 2$_*$, for the case $h _B + 1 = k_B$. Note that, in the case $h _B + 1 = k_B$, the above asymptotic expansion for function $C(\e)$ can not  be extended. Indeed, 
$\e^{h_B  - 2} C_2(\e)   = \e^{ h_B -2} ( C(\e) -  c_{h_C} \e^{- h_B} -  c_{h_C +1} \e^{- h_B +1})   
= - \frac{c_{h_C}}{b_{h_B} + b_{h_B + 1} \e + o_B(\e^{h_B+ 1}) \e^{- h_B}}$ $\times 
\frac{o_B(\e^{h_B +1}) \e^{- h_B -1} }{\e}$.  
The term  $\frac{o_B(\e^{h_B +1})\e^{- h_B -1}}{ \e}$ on the right hand side in the latter relation has  an uncertain asymptotic behavior as $\e \to 0$.

We can repeat the above arguments for the general case  $h _B + n = k_B$, for any $n = 0, 1, \ldots$ and to 
prove that, in the case $h _B + n = k_B$,  function $C(\e)$ can be represented in the form of $(h_C, k_C)$-expansion with parameters $h_C = - h_B, k_C = k_B - 2h_B = - h_B + n =  h_C + n$ and coefficients  
$c_{h_C}, \ldots, c_{k_C}$ given in proposition {\bf (iv)} of Lemma 2$_*$. Moreover,  
identity $B(\e) \cdot C(\e) \equiv 1, 0 < \e \leq \e'_0 $, let us find the corresponding remainder $o_C(\e^{k_C})$ from the following relation,
\begin{align}\label{nopty}
& (b_{h_B}\e^{h_B} + \cdots +  b_{k_B} \e^{k_B} + o_B(\e^{k_B})) \nonumber \\ 
& \times (c_{h_C}\e^{h_C}  + \cdots +  c_{h_C} \e^{k_C} + o_C(\e^{k_C})) \equiv 1.
\end{align}

Proposition {\bf (iii)} of Lemma 2$_*$,  applied to the product on the left hand side in relation (\ref{nopty}), permits to represent this product in the form of $(h, k)$-expansion 
with parameters $h =  h_B + h_C = h_B - h_B = 0$ and $k = (k_B +  h_C) \wedge (h_B + k_C)   = (k_B - h_B) \wedge (k_B - 2h_B + h_B) = k_B - h_B$.
By canceling coefficient for $\e^l$ on the left and right hand sides in  relation (\ref{nopty}), for $l = 0, \ldots, k_B - h_B$, and then, by solving equation  (\ref{nopty}) with respect to  the remainder $o_C(\e^{k_C})$, we get the formula for this remainder given in proposition {\bf (iv)} of Lemma 2$_*$, 
\begin{align}\label{foper}
o_C(\e^{k_C}) &  = - \frac{  \sum_{k_B - h_B < i + j, h_B \leq i \leq k_B, h_C \leq j \leq k_C} b_i c_j \e^{i + j} + \sum_{h_C \leq j \leq k_C}  c_j \e^j o_B(\e^{k_B}) }{b_{h_B}\e^{h_B} + \cdots +  b_{k_B} \e^{k_B} + o_B(\e^{k_B})}  \nonumber \\
& =  - \frac{  \sum_{k_B - h_B < i + j, h_B \leq i \leq k_B, h_C \leq j \leq k_C} b_i c_j \e^{i + j - h_B}}{b_{h_B} + \cdots +  b_{k_B} \e^{k_B - h_B} + o_B(\e^{k_B}) \e^{ h_B}} \nonumber \\
& \quad   - \frac{\sum_{h_C \leq j \leq k_C}  c_j \e^{j - h_B}o_B(\e^{k_B})}{b_{h_B} + \cdots +  b_{k_B} \e^{k_B - h_B} + o_B(\e^{k_B}) \e^{- h_B}}.
\end{align}

The assumptions made in proposition {\bf (iv)} of Lemma 2, imply that $B(\e) \neq 0$ and the following inequality holds for $0 < \e \leq \e_C$, where $\e_C $  is given  in proposition {\bf (iv)} of Lemma 2,
\begin{align}\label{nert}
& |b_{h_B} +  b_{h_B +1} \e + \cdots +  b_{k_B} \e^{k_B - h_B} + o_B(\e^{k_B}) \e^{- h_B}| \nonumber \\
&\quad \geq |b_{h_B}| - (|b_{h_B +1}| \e + \cdots +  |b_{k_B}| \e^{k_B - h_B} + G_B\e^{k_B - h_B + \delta_B}) 
 \nonumber \\
&\quad \geq \frac{|b_{h_B}|}{2} > 0.  
\end{align}

The existence of $\e'_0$ declared in proposition {\bf (iv)} of Lemma 2 is obvious. For example, one can choose $\e'_0 = \e_C $. It is also useful to note that formulas  given  in proposition {\bf (iv)} of Lemma 2 imply that  $\e_C = \e_B \wedge \tilde{\e}_B \in (0, \e_0]$, since $\e_B \in (0, \e_0]$ and $\tilde{\e}_B \in (0, \infty)$.    

The assumptions made in proposition {\bf (iv)} of Lemma 2 and inequality (\ref{nert})   imply that the following inequality holds, for $0 < \e \leq \e_C$, 
\begin{align}\label{bera}
|o_C(\e^{k_C})| & \leq  \e^{k_B - 2h_B + \delta_B}  (\frac{|b_{h_B}|}{2})^{-1}  \nonumber \\ 
& \quad  \times \big( \sum_{k_B - h_B < i + j, h_B \leq i \leq k_B, h_C \leq j \leq  k_C} |b_i| |c_j| \e_C^{i + j - k_B + h_B - \delta_B} \nonumber \\
& \quad \quad  + G_B \sum_{h_C \leq j \leq k_C} |c_j| \e_C^{j + h_B} \big).  
\end{align}

Inequality (\ref{bera}) proofs proposition {\bf (iv)} of Lemma 2.

Propositions {\bf (v)} of Lemmas 2$_*$ and 2 and relations {\bf (a)} -- {\bf (c)} given in these propositions  can be obtained by direct application, respectively, of propositions {\bf (iii)} and  {\bf (iv)} of  Lemmas 2$_*$ and 2, to the product $D(\e) = A(\e) \cdot \frac{1}{B(\e)}$.

Now, when it is already known that $D(\e) = A(\e) \cdot \frac{1}{B(\e)} $ is  a $(h_D, k_D)$-expansion,  with parameters $h_D = h_A - h_B$ and $k_D =  (k_A - h_B)  \wedge (k_B - 2 h_B + h_A)$, multiplication of $D(\e)$ by $B(\e)$ yields the following relation holding for $\e \in (0, \e'_0]$, 
\begin{align}\label{bertoma}
A(\e)  = D(\e) B(\e)  & = \, a_{h_A}\e^{h_A} + \cdots +  a_{h_A} \e^{k_A} + o_A(\e^{k_A})  \nonumber \\
& =  ( d_{h_D}\e^{h_D} + \cdots +  d_{h_D} \e^{k_D} + o_D(\e^{k_D})) \nonumber \\
& \quad \times  (b_{h_B}\e^{h_B} + \cdots +  b_{h_B} \e^{k_B} + o_B(\e^{k_B})). 
\end{align} 

By equating coefficients for powers $\e^l$ for $l = h_D, \ldots, k_D$ on the left and right hand sides of the third equality in  relation (\ref{bertoma}), we get  alternative formulas {\bf (e)} for coefficients $d_{h_d}, \ldots, d_{k_D}$ given in proposition {\bf (v)} of Lemma 2$_*$.

Proposition {\bf (iii)} of Lemma 2,  applied to the product on the right hand side in  (\ref{bertoma}), permits to represent this product in the form of $(h, k)$-expansion 
with parameters $h =  h_B + h_D = h_B + h_A - h_B = h_A$ and $k = (k_D + h_B) \wedge (k_B +  h_D) = ((k_A - h_B) \wedge (k_B - 2 h_B + h_A) + h_B) \wedge (k_B + h_A - h_B)  = k_A \wedge (k_B + h_A - h_B)$.
By canceling coefficient for $\e^l$ on the left and right hand sides in relation  (\ref{bertoma}), for $l = h_A, \ldots, k_A \wedge (k_B + h_A - h_B)$, and then, by solving equation  (\ref{bertoma}) with respect to  the remainder $o_D(\e^{k_D})$, we get the formula {\bf (f)} for this remainder given in proposition {\bf (v)} 
of Lemma 2$_*$,
\begin{equation*} 
\begin{aligned}
%\begin{align}\label{boputy}
o_D(\e^{k_D}) &  =  \frac{  \sum_{k_A \wedge (k_B + h_A - h_B) < l \leq k_A} a_l \e^{l} + o_A(\e^{k_A}) }{b_{h_B}\e^{h_B} + \cdots +  b_{k_B} \e^{k_B} + o_B(\e^{k_B})} \makebox[35mm]{}  \nonumber \\
\end{aligned}
\end{equation*}
\begin{align}\label{boputy}
&  \quad   - \frac{ \sum_{k_A \wedge (k_B + h_A - h_B) < i + j, h_B \leq i \leq k_B, h_D \leq j \leq k_D} b_i d_j \e^{i + j} }{b_{h_B}\e^{h_B} + \cdots +  b_{k_B} \e^{k_B} + o_B(\e^{k_B})}  \nonumber \\
& \quad   - \frac{ \sum_{h_D \leq j \leq k_D}  d_j \e^j o_B(\e^{k_B}) }{b_{h_B}\e^{h_B} + \cdots +  b_{k_B} \e^{k_B} + o_B(\e^{k_B})} \nonumber \\
& =  \frac{  \sum_{k_A \wedge (k_B + h_A - h_B) < l \leq k_A} a_l \e^{l - h_B} + o_A(\e^{k_A}) \e^{- h_B}}{b_{h_B} + \cdots +  b_{k_B} \e^{k_B - h_B} + o_B(\e^{k_B})\e^{- h_B} }  \nonumber \\
&    \quad   - \frac{ \sum_{k_A \wedge (k_B + h_A - h_B) < i + j, h_B \leq i \leq k_B, h_D \leq j \leq k_D} b_i d_j \e^{i + j - h_B} }{b_{h_B} + \cdots +  b_{k_B} \e^{k_B - h_B} + o_B(\e^{k_B})\e^{- h_B} }  \nonumber \\
&  \quad   - \frac{ \sum_{h_D \leq j \leq k_D}  d_j \e^{j - h_B}o_B(\e^{k_B})  }{b_{h_B} + \cdots +  b_{k_B} \e^{k_B - h_B} + o_B(\e^{k_B})\e^{- h_B}}.
\end{align}

Inequality (\ref{nert}) and the assumptions made in proposition {\bf (v)} of Lemma 2 finally  imply that the following inequality holds, for $0 < \e \leq \e_D$  given in relation {\bf (f)} of this proposition, 
\begin{align}\label{dwar}
|o_D(\e^{k_D})| & \leq  \e^{k_D  + \delta_D}  (\frac{|b_{h_B}|}{2})^{-1}   \big( \sum_{k_A \wedge (k_B + h_A - h_B) < l \leq  k_A} |a_l| \e_D^{l - k_D   - h_B - \delta_D} \nonumber \\
& \quad   +  \sum_{k_A \wedge (k_B + h_A - h_B) < i + j, h_B \leq i \leq k_B, h_D \leq j \leq  k_D} |b_i| |d_j| \e_D^{i + j - k_D - h_B - \delta_D} \nonumber \\
& \quad   + G_A e_D^{k_A + \delta_A - h_B - k_D  - \delta_D}  \nonumber \\
& \quad  +  G_B \sum_{h_D \leq j \leq k_D} |d_j| \e_D^{j +k_B + \delta_B  - h_B - k_D  - \delta_D} \big).  
\end{align}

Inequality  (\ref{dwar})  yields relations {\bf (d)} -- {\bf (f)} given in proposition {\bf (v)} of Lemma 2.  

{\bf A.3}. Lemma 3$_*$ is a direct corollary of Lemma 2$_*$. Proofs of propositions {\bf (i)} and  {\bf (ii)} in Lemma 3 are analogous to  proofs of  the corresponding propositions in Lemma  2. Proposition {\bf (iii)} of Lemma 3 is obvious.

{\bf A.4}. The first two identities for Laurent asymptotic expansions given in proposition {\bf (i)} of Lemma 4$_*$ are obvious. The third identity given in 
this proposition follows in an obvious way from proposition {\bf (i)} of Lemma 2$_*$. By applying propositions {\bf (iii)} and  
{\bf (iv)} of Lemma 2$_*$ to the product $C(\e) = A(\e) \cdot A(\e)^{-1}$,  we get  parameters $h_C = h_{A \cdot A^{-1}}  = h_A - h_A = 0, \, k_C = k_{A \cdot A^{-1}} = (k_A - h_A) \wedge (k_A - 2h_A + h_A) = k_A - h_A$ and  coefficients $c_n = {\rm I}(n = 0), n = 0, \ldots, k_C$.  Also, relations   (\ref{nopty}) and (\ref{foper})  imply that the elimination identity $A(\e) \cdot A(\e)^{-1} \equiv 1$ holds, since the remainder of Laurent asymptotic expansion for function $A(\e)^{-1}$  is given by formula {\bf (c)} from proposition {\bf (iv)} of Lemma 2$_*$. 

Propositions {\bf (ii)} and  {\bf (iii)} of Lemma 4$_*$ in the parts concerned commutative property of summation and multiplication operations  follow  from, respectively,  propositions  {\bf (ii)} and {\bf (iii)}  of Lemma  2$_*$.

Let  $D(\e) = (A(\e) + B(\e)) +  C(\e) = A(\e) +  (B(\e) +   C(\e))$. Using propositions {\bf (ii)} of Lemma  2$_*$, we get,
$h_D  = h_{(A+ B) + C} = (h_A \wedge h_B) \wedge h_C = h_A \wedge (h_B \wedge h_C) = h_{A+ (B + C)}$ and
$k_D  = k_{(A+ B) + C} = (k_A \wedge k_B) \wedge k_C = k_A \wedge (k_B \wedge k_C) = k_{A+ (B + C)}$. These relations and Lemma 1$_*$ imply equalities for the corresponding coefficients and remainders,  for the asymptotic expansions of functions $(A(\e) + B(\e)) + C(\e)$ and $A(\e) +  (B(\e)  + C(\e))$. The above remarks prove proposition {\bf (ii)} of Lemma 4$_*$ in the part concerned with the associative property of summation operation for Laurent asymptotic expansions.

Let $D(\e) = (A(\e) \cdot B(\e)) \cdot  C(\e) = A(\e) \cdot  (B(\e) \cdot   C(\e))$. Using propositions {\bf (iii)} of Lemma  2$_*$, we get,
$h_D  = h_{(A \cdot  B) \cdot  C} = h_{A \cdot  B}  +  h_C = h_A + h_B + h_C = h_A + h_{B \cdot C} = h_{A \cdot  (B \cdot  C)}$ and 
$k_D =   k_{(A \cdot  B) \cdot  C} =  (k_{A \cdot B} + h_C) \wedge (k_C + h_{A \cdot B}) = ((k_A + h_B) \wedge (k_B + h_A)) + h_C) \wedge (k_C + (h_A + h_B)) = 
(k_A + h_B + h_C) \wedge (k_B + h_A + h_C) \wedge (k_C + h_A + h_B) = (k_A + (h_B + h_C)) \wedge ((k_B + h_C) \wedge (k_C + h_B))  + h_A) =
(k_A + h_{B \cdot C}) \wedge (k_{B \cdot C} + h_A) = k_{A \cdot ( B \cdot  C)}$.  These relations and Lemma 1$_*$ imply equalities  for the corresponding coefficients and remainders, for the asymptotic expansions of functions $(A(\e) \cdot B(\e)) \cdot C(\e)$ and $A(\e) \cdot (B(\e) \cdot C(\e))$. 
The above remarks prove proposition {\bf (iii)} of Lemma 4$_*$ in the part concerned with the associative property of multiplication operation for Laurent asymptotic expansions.

Let  $D(\e)  = (A(\e) + B(\e)) \cdot C(\e) = A(\e)\cdot C(\e) + B(\e) \cdot C(\e)$.  Using propositions {\bf (ii)} and  {\bf (iii)} of Lemma  2$_*$, we get,  $h_D = h_{(A + B) \cdot C} =  h_{A + B} + h_C = h_A \wedge h_B + h_C = (h_A + h_C) \wedge (h_B + h_C) = h_{A \cdot C} \wedge h_{B \cdot C}  = 
h_{A\cdot C + B \cdot C}$ and   $k_D = k_{(A + B) \cdot C} =  (k_{A + B} + h_C) \wedge (k_C + h_{A + B}) =   (k_A \wedge k_B + h_C) \wedge (k_C + h_A \wedge h_B) = (k_A + h_C) \wedge (k_B + h_C) \wedge (k_C + h_A) \wedge (k_C + h_B) = ((k_A + h_C) \wedge (k_C + h_A)) \wedge ((k_B + h_C) \wedge  (k_C + h_B)) = 
k_{A \cdot C} \wedge k_{B \cdot C} = k_{A \cdot C + B \cdot C}$.  
These relations and Lemma 1$_*$ imply equalities for the corresponding coefficients and remainders, for the asymptotic expansions of functions $(A(\e) + B(\e)) \cdot C(\e)$ and $A(\e)\cdot C(\e) + B(\e) \cdot C(\e)$. The above remarks prove proposition {\bf (iv)} of Lemma 4$_*$ concerned with the distributive property of summation and  multiplication operations for Laurent asymptotic expansions.

Propositions {\bf (i)} and  {\bf (ii)} of Lemma 4 readily follow from, respectively,  
propositions  {\bf (ii)} and {\bf (iii)}  of Lemma  2. Finally, proposition {\bf (iii)} of Lemma 4 follows from relations 
$\delta_{cA}= \delta_A$, $\delta_{A + B}, \delta_{A \cdot B}$, $\delta_{A / B} \geq \delta_A \wedge \delta_B$ and $\delta_{A_1 + \cdots + A_N}$,  
$\delta_{A_1 \times \cdots \times A_N} \geq \min_{1 \leq m \leq N} \delta_{A_m}$, given, respectively, in Lemmas 2 and 3. $\Box$ \\

{\bf Appendix B: Examples} \\ 

Let us, first,  comment some general questions connected with  construction of examples illustrating the asymptotic results presented in the paper.

Let $\YY_i \neq \emptyset, i \in \XX$ be some subsets of space $\XX$  such that condition ${\bf A}$ {\bf (c)} holds for these sets, i.e.,   for every pair of states $i, j \in \XX$, there exists an integer  $n_{ij} \geq 1$ and a chain of states $i = l_{ij, 0}, l_{ij, 1}, \ldots, l_{ij, n_{ij}} = j$ such that  $l_{ij, 1} \in \YY_{l_{ij, 0}}, \ldots, l_{ij, n_{ij}} \in  \YY_{l_{ij, n_{ij}-1}}$.

Let us also choose some $\e_0 \in (0, 1]$. 

We define $p_{ij}(\e)  = 0,  \e \in (0, \e_0], j \in  \overline{\YY}_i, \, i \in \XX$, i.e., assume that condition  
${\bf A}$ {\bf (b)}  holds.

Let $p_{ij}(\e), \e \in (0, \e_0], j \in  \YY_i, 
i \in \XX$ be some real-valued functions which satisfy condition ${\bf D}$, i.e., can be represented in the form of Taylor asymptotic  expansions, $p_{ij}(\e) =  \sum_{l = l_{ij}^-}^{l_{ij}^+} a_{ij}[l]\e^l + o_{ij}(\e^{l_{ij}^+}), \, \e \in (0, \e_0]$, for $j \in \YY_i, i \in \XX$, where {\bf (a)}  $a_{ij}[ l_{ij}^-] > 0$ and $0 \leq l_{ij}^- \leq l_{ij}^+ < \infty$, for $j \in \YY_i, i \in \XX$;  {\bf (b)}  
$o_{ij}(\e^{l_{ij}^+}) /\e^{l_{ij}^+} \to 0$ as $\e \to 0$, for $j \in \YY_i, i \in \XX$.  

Condition ${\bf D}$ does not guarantee that matrix  $\| p_{ij}(\e) \|$ is stochastic, for every  $\e \in (0, \e_0]$. This can be achieved by imposing some additional 
conditions on coefficients and remainders in the above asymptotic expansions. 

Let us choose arbitrary  numbers $0 < \alpha_{ij} < \frac{1}{2}$, for $j \in \YY_i, i \in \XX$.

Condition ${\bf D}$ {\bf (b)} guarantees that, for every $j \in  \YY_i, \, 
i \in \XX$,    there exists $\e_{\alpha_{ij}, ij} \in (0, \e_0]$ such that, for $\e \in (0, \e_{\alpha_{ij}, ij}]$,
\begin{equation}\label{condas}
| o_{ij}(\e^{l_{ij}^+}) / \e^{l_{ij}^+} | \leq \alpha_{ij} | a_{ij}[l_{ij}^-] |.
\end{equation}

It is useful to note that in the case, where condition ${\bf D'}$ holds, an explicit value for parameters 
$\e_{\alpha_{ij}, ij}$ can be derived from the inequalities,   $|o_{ij}(\e^{l_{ij}^+})| \leq G_{ij}\e^{l_{ij}^+ + \delta_{ij}}, \e \in (0, \e_{ij}], j \in  \YY_i, \, i \in \XX$,  penetrating this condition. 
These inequalities imply  that relation (\ref{condas}) holds for 
$\e_{\alpha, ij} = \e_{ij} \wedge (\frac{\alpha_{ij}  | a_{ij}[l_{ij}^-] |}{G_{ij}})^{\frac{1}{\delta_{ij}}}$, $j \in  \YY_i, \, 
i \in \XX$.  

Let us also define $A_{\e_0, ij} = \sum_{l_{ij}^- < l \leq l_{ij}^+} |a_{ij}[l] | \e_0^{l - l_{ij}^-  -1}, j \in  \YY_i, \, 
i \in \XX$ and $\e'_{\alpha_{ij},  ij} =  \e_{\alpha_{ij}, ij}$ if $A_{\e_0, ij} = 0$ or 
$\e'_{\alpha_{ij},  ij} = \e_{\alpha_{ij}, ij} \wedge \frac{\alpha_{ij} | a_{ij}[l_{ij}^- ] |}{A_{\e_0, ij}}$ 
if $A_{\e_0, ij} > 0$.
 
The following inequality holds, for every $\e \in (0, \e'_{\alpha_{ij},  ij}]$ and $j \in  \YY_i, \, 
i \in \XX$,
\begin{align}\label{pose}
p_{ij}(\e) & \geq   \e^{l_{ij}^-}  (a_{ij}[l_{ij}^- ] - \e A_{\e_0, ij} - \e_0^{l_{ij}^+ - l_{ij}^-} 
| o_{ij}(\e^{l_{ij}^+}) /\e^{l_{ij}^+} |) \nonumber \\
& \geq     \e^{l_{ij}^-} a_{ij}[l_{ij}^- ] (1 - 2 \alpha_{ij}) > 0. 
\end{align}

Let us now define $\e'_0 = \min_{j \in  \YY_i, \, i \in \XX} \e'_{\alpha_{ij},  ij}$. 

Obviously, $p_{ij}(\e) > 0, j \in  \YY_i, \, i \in \XX, \e \in (0, \e'_0]$. Thus, conditions  ${\bf A}$ {\bf (a)}  -- {\bf (b)} hold, if  parameter $\e_0$ is replaced by the new value $\e'_0$.

The question about holding of the stochasticity relation $ \sum_{j \in  \YY_{i}}  p_{ij}(\e) = 1, \, \e \in (0, \e_0], i \in \XX$ is 
more complex.  According Lemma 5$_*$, under conditions   ${\bf A}$ {\bf (a)}  -- {\bf (b)} and ${\bf D}$, the above stochasticity relation is equivalent to condition ${\bf F}$ formulated in Section 3$_*$.

First, condition ${\bf F}$ requires holding of the following relation,   
\begin{equation}\label{retaw}
 \sum_{j \in \YY_i } a_{ij}[l] = {\rm I}(l = 0), \
0  \leq l \leq l_{i, \YY_i}^+, \ i \in \XX, 
\end{equation}
where (a) $ l_{i, \YY_i}^{\pm} = \min_{j \in \YY_i} l_{ij}^{\pm}, i \in \XX$ and (b)   
$a_{ij}[l] = 0$, for $0 \leq l < l_{ij}^-, j \in \YY_i, i \in \XX$.

Note that relation (\ref{retaw}) implies that parameters $l_{i, \YY_i}^{-}  = 0, i \in \XX$.

It is not difficult to choose coefficients $a_{ij}[l], l = l_{ij}^- \leq l \leq l_{ij}^+, j \in \YY_i, \ i \in \XX$ in such way 
that relation  (\ref{retaw}) would hold. Any such coefficients, with the first coefficients $a_{ij}[ l_{ij}^-] > 0, j \in \YY_i, \ 
i \in \XX$,  can serve as coefficients in the asymptotic expansions penetrating condition ${\bf D}$.

Second, condition ${\bf F}$ requires holding of the following identity, for every  $i \in \XX$,
\begin{equation}\label{hopter}
\sum_{j \in \YY_i} ( \sum_{l_{i, \YY_i}^+ 
<  l \leq l_{ij}^+} a_{ij}[l]\e^l + o_{ij}(\e^{l_{ij}^+}) ) \equiv 0.
\end{equation}

Remainders $o_{ij}(\e^{l_{ij}^+}), j \in \YY_i, \ i \in \XX$ satisfying the above identities can be chosen in different ways.

The simplest one is to choose $o_{ij}(\e^{l_{ij}^+}) \equiv 0, \, j \in \YY_i, \,  i \in \XX$. In this case, the above identities would reduce to equalities,  $\sum_{j \in \YY_i}  a_{ij}[l]  = 0, l_{i, \YY_i}^+ <  l \leq l^*_{i, \YY_i}, i \in \XX$, where  (a) $ l^*_{i, \YY_i} = \max_{j \in \YY_i} l_{ij}^{+}, i \in \XX$ and (b)   $a_{ij}[l] = 0$, for $l_{ij}^+ < l \leq  l^*_{i,  \YY_i}, j \in \YY_i, i \in \XX$.
These equalities supplement equalities given in relation (\ref{retaw}). Such choice of remainders corresponds to models with polynomial perturbations. 

We, however, would like to impose on remainders  conditions mainly required of them by conditions ${\bf D}$ or ${\bf D'}$.

There always exist $j_i \in \YY_i, i \in \XX$ such that $l^+_{ij_i} =  l_{i, \YY_i}^{+}, i \in \XX$.

Identity (\ref{hopter}) can be rewritten in the following form, for every  $i \in \XX$,
\begin{equation}\label{hoptera}
o_{ij_i}(\e^{l_{ij_i}^+}) \equiv - \sum_{j \in \YY_i , j \neq j_i} ( \sum_{l_{i, \YY_i}^+ 
<  l \leq l_{ij}^+} a_{ij}[l]\e^l + o_{ij}(\e^{l_{ij}^+}) ).
\end{equation}

Relation (\ref{hoptera}) can be used as the formula defining remainders $o_{ij_i}(\e^{l_{ij_i}^+}), i \in \XX$,  via remainders 
$o_{ij}(\e^{l_{ij}^+}), j \in \YY_i , j \neq j_i, i \in \XX$ penetrating the corresponding asymptotic expansions 
in condition ${\bf D}$.

Since $l^+_{ij_i} =  l_{i, \YY_i}^{+}, i \in \XX$,  the following relation holds, for remainder $o_{ij}(\e^{l_{ij}^+})$ defined by relation  (\ref{hoptera}),  for every $i \in \XX$,
\begin{equation}\label{hopteramo}
o_{ij_i}(\e^{l_{ij_i}^+}) /\e^{l_{ij_i}^+} \to 0 \ {\rm as} \ \e \to 0.
\end{equation}

Thus, remainders $o_{ij_i}(\e^{l_{ij_i}^+}), i \in \XX$ defined by relation (\ref{hoptera}) can also serve in the corresponding asymptotic expansions in condition ${\bf D}$.

Moreover, let us assume that remainders $o_{ij}(\e^{l_{ij}^+}), j \in \YY_i , j \neq j_i, i \in \XX$  satisfy the inequalities, $|o_{ij}(\e^{l_{ij}^+})| \leq G_{ij}\e^{l_{ij}^+ + \delta_{ij}}, \e \in (0, \e_{ij}], j \in  \YY_i, j \neq i_j, \, i \in \XX$,  penetrating condition ${\bf D'}$.

Let us define $\e_{ij_i} = \min_{j \in \YY_i, j \neq j_i} \e_{ij}, i \in \XX$ and  $\delta_{ij_i} = 
\min_{j \in \YY_i, j \neq j_i} \delta_{ij}, i \in \XX$.

In this case, the following inequality holds, for  every $\e \in (0, \e_{ij_i}], i \in \XX$,
\begin{align}\label{hopteraba}
| o_{ij_i}(\e^{l_{ij_i}^+}) | & \leq   \sum_{j \in \YY_i , j \neq j_i} ( \sum_{l_{i, \YY_i}^+ 
<  l \leq l_{ij}^+} |a_{ij}[l] | \e^l  +  |o_{ij}(\e^{l_{ij}^+})| ) \nonumber \\
& \leq  \big( \sum_{j \in \YY_i , j \neq j_i}  ( \sum_{l_{i, \YY_i}^+ 
<  l \leq l_{ij}^+} |a_{ij}[l] | \e_0^{l - l_{ij_i}- \delta_{ij_i}}    \nonumber \\ 
& \quad +   \e_0^{l_{ij}^+ - l_{ij_i}^+ + \delta_{ij} -\delta_{ij_i}}  G_{ij}) \big )  \e^{l_{ij_i}^+  + \delta_{ij_i}}
=  G_{ij_i}  \e^{l_{ij_i}^+  + \delta_{ij_i}}.  
\end{align}

Thus, the inequalities,  $|o_{ij_i}(\e^{l_{ij_i}^+})| \leq G_{ij_i}\e^{l_{ij_i}^+ + \delta_{ij_i}}, \e \in (0, \e_{ij_i}], i \in \XX$, penetrating condition ${\bf D'}$ hold for remainders $o_{ij_i}(\e^{l_{ij_i}^+}), i \in \XX$, with parameters  $ \e_{ij_i}, \delta_{ij_i}$ and 
$G_{ij_i}$ defined above.

As follows from the above remarks, identity (\ref{hoptera}) holds for remainders  $o_{ij}(\e^{l_{ij}^+}), j \in \YY_i, , i \in \XX$ , for $\e \in (0, \e'_0], i \in \XX$, where  $\e'_0 =  \min_{j \in \YY_i} \e_{ij} = \min_{j \in \YY_i, j \neq j_i} \e_{ij}$. Thus, condition  ${\bf F}$ holds, if  parameter $\e_0$ is replaced by the new value $\e'_0$. In this case, functions $p_{ij}(\e), i, j \in \XX$ can, for every $\e \in (0, \e'_0]$, serve as transition probabilities of a Markov chain.
 
Note that remainders $o_{ij}(\e^{l_{ij}^+}), j \in \YY_i, i \in \XX$ constructed above can be very irregular functions.
Let us, for example,  consider the case,  where all asymptotic expansions in condition ${\bf D}$ have the same order, i.e., 
parameters $l^+_{ij} = l^+, j \in \YY_i, i \in \XX$. In this case, identities  (\ref{hopter}) take the form, 
 $\sum_{j \in \YY_i} o_{ij}(\e^{l^+}) = 0, \e \in (0, \e_0], i \in \XX$. Condition ${\bf D}$ requires that $o_{ij}(\e^{l^+})/ \e^{l^+} \to 0$ as $\e \to 0$, for $j \in \YY_i, i \in \XX$. Remainders $o_{ij}(\e^{l^+}), j \in \YY_i, i \in \XX$  can be  continuous functions of $\e$ taking zero value in at most finite numebrs of points. However, let us multiply them, for example, by the Dirichlet function $D(\e)$.  The new remainders   $o'_{ij}(\e^{l^+}) = D(\e)o_{ij}(\e^{l^+}), j \in \YY_i, i \in \XX$ also satisfy identities (\ref{hopter}) and $o'_{ij}(\e^{l^+})/ \e^{l^+} \to 0$ as $\e \to 0$, for $j \in \YY_i, i \in \XX$. At the same time, they are  very irregular functions. This example is, of course, an artificial one. But, it well illustrates the above statement about possible irregularity of remainders and, in sequel,  transition probabilities,  as functions of the perturbation parameter.

Let us also make some remarks concerned the expected sojourn times.

First, let us define $e_{ij}(\e)  = 0, \e \in (0, \e_0] \, j \in  \overline{\YY}_i, \, i \in \XX$ that is consistent with  condition  
${\bf A}$ {\bf (b)}.

Let us also $e_{ij}(\e), \e \in (0, \e_0] \, j \in  \YY_i, \, i \in \XX$ be some real-valued functions which satisfy condition ${\bf E}$, i.e., can be represented in the form of  Laurent  asymptotic expansions, $e_{ij}(\e) =  \sum_{l = m_{ij}^-}^{m_{ij}^+} b_{ij}[l]\e^l + \dot{o}_{ij}(\e^{m_{ij}^+}), \, \e \in (0, \e_0]$, for $j \in \YY_i, \,  i \in \XX$, where {\bf (a)} $b_{ij}[ m_{ij}^-] > 0$ and $-\infty < m_{ij}^- \leq m_{ij}^+ < \infty$, for $j \in \YY_i, i \in \XX$; {\bf (b)} $ \dot{o}_{ij}(\e^{m_{ij}^+})/\e^{m_{ij}^+} \to 0$ as $\e \to 0$, for 
$j \in \YY_i, i \in \XX$. 

Condition ${\bf E}$ {\bf (b)} guarantees that, for every $j \in  \YY_i, \, 
i \in \XX$   there exists $\dot{\e}_{\alpha_{ij}, ij} \in (0, \e_0]$ such that, for $\e \in (0, \dot{\e}_{\alpha_{ij}, ij}]$,
\begin{equation}\label{condasbes}
| \dot{o}_{ij}(\e^{m_{ij}^+}) / \e^{m_{ij}^+} | \leq \alpha_{ij} | b_{ij}[m_{ij}^-] |.
\end{equation}

It is useful to note that in the case, where condition ${\bf E'}$ is assumed to hold, an explicit value for parameters 
$\dot{\e}_{\alpha_{ij}, ij}$ can be derived from the inequalities,   $|\dot{o}_{ij}(\e^{m_{ij}^+})| \leq \dot{G}_{ij}\e^{m_{ij}^+ + \dot{\delta}_{ij}}, \e \in (0, \e_{ij}], j \in  \YY_i, \, i \in \XX$,  penetrating this condition. 
These inequalities yield  that relation (\ref{condasbes}) holds for 
$\dot{\e}_{\alpha, ij} = \dot{\e}_{ij} \wedge (\frac{\alpha_{ij}  | b_{ij}[m_{ij}^-] |}{\dot{G}_{ij}})^{\frac{1}{\dot{\delta}_{ij}}}, j \in  \YY_i, \,  i \in \XX$.  

Let us also define $B_{\e_0, ij} = \sum_{m_{ij}^- < l \leq m_{ij}^+} |b_{ij}[l] | \e_0^{l - m_{ij}^-  -1}, j \in  \YY_i, \, 
i \in \XX$ and $\dot{\e}''_{\alpha_{ij},  ij} =  \dot{\e}_{\alpha_{ij}, ij}$ if $B_{\e_0, ij} = 0$ or 
$\dot{\e}''_{\alpha_{ij},  ij} = \dot{\e}_{\alpha_{ij}, ij} \wedge \frac{\alpha_{ij} | b_{ij}[l_{ij}^- ] |}{B_{\e_0, ij}}$ 
if $B_{\e_0, ij} > 0$.
 
The following inequality holds, for every $\e \in (0, \dot{\e}''_{\alpha_{ij}, ij}]$ and $j \in  \YY_i, \, 
i \in \XX$,
\begin{align}\label{poseva}
e_{ij}(\e) & \geq   \e^{m_{ij}^-}  (b_{ij}[m_{ij}^- ] - \e B_{\e_0, ij} - \e_0^{m_{ij}^+ - m_{ij}^-} 
| o_{ij}(\e^{m_{ij}^+}) /\e^{m_{ij}^+} |) \nonumber \\
& \geq     \e^{m_{ij}^-}  b_{ij}[m_{ij}^- ] (1 - 2 \alpha_{ij}) > 0. 
\end{align}

Let us now define $\e''_0 = \min_{j \in  \YY_i, \, i \in \XX} \dot{\e}''_{\alpha_{ij},  ij}$. 

Obviously, $e_{ij}(\e) > 0, j \in  \YY_i, \, i \in \XX, \e \in (0, \e''_0]$. This is consistent, with conditions ${\bf A}$ {\bf (a)} and ${\bf B}$. 

Finally, let us define $\tilde{\e}_0 = \e'_0 \wedge  \e''_0$. Parameter, $\tilde{\e}_0$  can serve as a new value for parameter $\e_0$.

Functions  $p_{ij}(\e), i, j \in \XX$ and  $e_{ij}(\e), i, j \in \XX$ constructed above  can serve, respectively, as  transition probabilities of the embedded Markov chain $\eta^{(\e)}_n$  and expectations of sojourn times for some semi-Markov process $\eta^{(\e)}(t)$, for every $\e \in (0, \tilde{\e}_0]$.  A variant of transition probabilities for such  semi-Markov processes is given in Section 3$_*$. 

In conclusion, let us  consider a numerical example. 

We assume that the semi-Markov process $\eta^{(\e)}(t)$ is, for every $\e \in (0, \e_0]$, a semi-Markov process 
with the phase space $\XX = \{1, 2, 3 \}$. 

The transition sets are $\YY_1 = \{1, 2 \}, \YY_2 = \{1, 2, 3 \}, \YY_3 = \{1, 2\}$. 

The  $3 \times 3$ matrix  of transition probabilities $\| p_{ij}(\e) \|$, for the corresponding embedded Markov chain $\eta^{(\e)}_n$, 
has the following form, 
{\small
\begin{equation}\label{partysak}
\left \| 
\begin{array}{llllll}
1 - \e^2 - \e^3 + o_{11}(\e^3) & \e^2 + \e^3 + o_{12}(\e^3)  & 0  \vspace{1mm} \\
\frac{\e}{2} + \frac{\e ^2}{2} - \e^3 + o_{21}(\e^3)    & 1 - \e - \e^2  + o_{22}(\e^2)   &  \frac{\e}{2} 
+ \frac{\e ^2}{2} + 2 \e^3 + o_{23}(\e^3)  \vspace{1mm}  \\
\frac{1}{2} + \e^2 - \e^3 + o_{31}(\e^3)  & \frac{1}{2} - \e^2 + \e^3 + o_{32}(\e^3) & 0 \\
\end{array}
\right \|. 
\end{equation}
}

The $3 \times 3$ matrix of expectations of sojourn times $\| e_{ij}(\e) \|$, for the semi-Markov process $\eta^{(\e)}(t)$,  has the 
following form, 
{\small
\begin{equation}\label{partys}
\left \| 
\begin{array}{lllllll}
\e + \e^2 + \dot{o}_{11}(\e^2) & \e^3 + \e^4 + \dot{o}_{12}(\e^4)  & 0 \vspace{1mm}  \\
\e + \e ^2  + \dot{o}_{21}(\e^2)    & 1 - \e + \e^3  + \dot{o}_{22}(\e^3)   & 2 + \e +  
\e ^2  + \dot{o}_{23}(\e^2)  \vspace{1mm}  \\
\e^{-1} + 1 + \dot{o}_{31}(1)  & 2 \e^{-1} + \e + \dot{o}_{32}(\e) & 0 \\
\end{array}
\right \|.
\end{equation}
}

In the asymptotic expansions penetrating relations (\ref{partysak}) and (\ref{partys}), the coefficients $a_{ij}[l_{ij}^-], b_{ij}[m_{ij}^-] > 0, \, j \in \YY_i, \ i \in \XX$, and  coefficients  $a_{ij}[l],  l = l_{ij}^-, \ldots,  l_{ij}^+, j \in \YY_i, \ i \in \XX$ satisfy relation  (\ref{retaw}). We also assume that parameter $\e_0 = \tilde{\e}_0$  and  remainders $o(\e^{l_{ij}^+}), \dot{o}(\e^{m_{ij}^+}),  j \in \YY_i, \ i \in \XX$, in the asymptotic expansions representing elements of  matrices given in relations (\ref{partysak}) and (\ref{partys}),  are chosen according the procedures described above, in particular, the identities (\ref{hopter}) hold. In this case, matrices, given in the above relations,  can, for every $\e \in (0, \e_0]$,  serve as, respectively,  the matrix of transition probabilities for the corresponding embedded Markov chain and the matrix of expectations of sojourn times,  for the semi-Markov process $\eta^{(\e)}(t)$, and  conditions ${\bf A}$ --  ${\bf E}$ hold.  

The matrices  of transition probabilities $\| p_{ij}(\e) \|$, for the embedded Markov 
chains $\eta^{(\e)}_n$, and $\| p_{ij}(0) \| $,  for the limiting  Markov 
chain $\eta^{(0)}_n$, have, respectively, the following forms,
{\small
\begin{equation}\label{partysa}
\left \| 
\begin{array}{lllllll}
\bullet & \bullet  & 0  \\
\bullet  & \bullet  & \bullet    \\
\bullet  &  \bullet  & 0 \\
\end{array}
\right \| \quad {\rm and} \quad   
\left \| 
\begin{array}{lllllll}
1 & 0 & 0  \\
0 & 1 & 0   \\
 \frac{1}{2} &  \frac{1}{2} & 0 \\
\end{array}
\right \|,
\end{equation}
}
where symbol $\bullet$ indicates positions of positive elements in  matrices $\| p_{ij}(\e) \|$, $\e \in (0, \e_0]$. 

The phase space $\XX$ is one class of communicative states for the 
embedded Markov chain $\eta^{(\e)}_n$, for every $\e \in (0, \e_0]$, while it consists of two closed classes of 
communicative states $\XX_1 = \{ 1 \},  \XX_2 = \{ 2 \}$ and the class of transient states $ \XX_3 = \{ 3 \}$,  for the 
limiting  Markov chain $\eta^{(0)}_n$. 

By excluding the state $1$ from the phase space $\XX$ and using the algorithm described in Section 5$_*$, we construct the reduced semi-Markov processes  $_1\eta^{(\e)}(t)$, with the phase space $_1\XX = \{ 2, 3 \}$. Conditions ${\bf A}$ --  ${\bf E}$  hold for these reduced semi-Markov processes. The corresponding transition sets are  $_1\YY_2 = \{ 2, 3 \}$ and $_1\YY_3 = \{ 2 \}$. By applying the algorithms described in Lemma 8$_*$ and Theorems 2$_*$ and 3$_*$, we  can compute the  $2 \times 2$ matrices $_1{\mathbf P} = \| _1p_{ij}(\e) \|$ and $_1{\mathbf E} =  \| _1e_{ij}(\e) \|$. These matrices  take the following forms, 
{\small
\begin{equation}\label{party}
_1{\mathbf P} = \left \| 
%\begin{array}{llllll}
%_1p_{22}(\e) & _1p_{23}(\e) \vspace{1mm}  \\
%_1p_{32}(\e) & _1p_{33}(\e) \\
%\end{array}
%\right \|  =
%\left \| 
\begin{array}{llllll}
1 - \frac{1}{2}\e - \frac{1}{2}\e^2 + \, _1o_{22}(\e^2) & \frac{\e}{2} + \frac{\e^2 }{2} +  2\e^3  + \, _1o_{23}(\e^3)  \vspace{1mm}  \\
1    & 0      \\
\end{array}
\right \|,  
\end{equation}
}
and 
{\small
\begin{equation}\label{party}
_1{\mathbf E}  = \left \| 
%\begin{array}{llllll}
%_1e_{22}(\e) & _1e_{23}(\e) \\
%_1e_{32}(\e) & _1e_{33}(\e) \\
%\end{array}
%\right \|  =
%\left \| 
\begin{array}{llllll}
\frac{3}{2} + \frac{1}{2}\e  + \, _1\dot{o}_{22}(\e) & 2 + \e + \e^2  + \, _1\dot{o}_{23}(\e^2)  \vspace{1mm}   \\
\frac{7}{2}\e^{-1} + 1 + \, _1\dot{o}_{32}(1)  &  0      \\
\end{array}
\right \|. 
\end{equation}
}

By excluding the state $2$ from the reduced phase space $_1\XX = \{ 2, 3 \}$,  we construct the ``final''  reduced semi-Markov processes  
$_{\langle1, 2 \rangle}\eta^{(\e)}(t)$, with the one-state phase space $_{\langle1, 2 \rangle}\XX = \{ 3 \}$. Conditions ${\bf A}$ --  ${\bf E}$ also hold for these semi-Markov processes. The corresponding transition set $_{\langle1, 2 \rangle}\YY = \{ 3 \}$. The transition probability 
$_{\langle1, 2 \rangle}p_{33}(\e) \equiv 1$. By applying the algorithms described in  Theorem 4$_*$, we  can also compute the Laurent asymptotic expansion for the expected return time,  $_{\langle1, 2 \rangle}e_{33}(\e) = E_{33}(\e) = \frac{21}{2}\e^{-1} - 3 + 
 \ddot{o}_{33}(1)$. The Laurent asymptotic expansion for the expected sojourn time $e_{3}(\e)  = e_{31}(\e) + e_{32}(\e) +  e_{33}(\e)$, obtained  with the use of the multiple summation rule given in Lemma 3$_*$, has the following  form,  $e_{3}(\e) = 3\e^{-1} + 1 + 
 \dot{o}_3(1)$. Finally, the algorithm described in Theorem 5$_*$ gives the following asymptotic expansion, for the stationary probability $\pi_3(\e) = \frac{e_{3}(\e)}{E_{33}(\e)} = \frac{3\e^{-1} + 1 + \dot{o}_3(1)}{\frac{21}{2}\e^{-1} - 3 +  \ddot{o}_{33}(1)} = \frac{2}{7} + \frac{26}{147}\e + o_3(\e)$.

Also, by excluding the state $3$ from the reduced phase space  $_1\XX = \{ 2, 3 \}$ and applying the algorithms described in Theorem 4$_*$, we  can compute the Laurent asymptotic expansion for the expected return time $_{\langle1, 3 \rangle}e_{22}(\e) = E_{22}(\e) = \frac{21}{4} + \frac{15}{4} \e  +  \ddot{o}_{22}(\e)$. In this case,  the asymptotic expansion for the expected sojourn time  $e_{2}(\e)  = 3 + \e + 2 \e^2 + \dot{o}_2(\e^2)$, and the algorithm described in Theorem 5$_*$ gives the following asymptotic expansion for the stationary probability $\pi_2(\e) = \frac{e_{2}(\e)}{E_{22}(\e)} = \frac{3 + \e + 2 \e^2 + \dot{o}_2(\e^2)}{ \frac{21}{4} + \frac{15}{4} \e  +  \ddot{o}_{22}(\e)} = \frac{4}{7} - \frac{32}{147}\e + o_2(\e)$.

As far as the stationary probability  $\pi_1(\e)$ is concerned, the corresponding asymptotic expansion can be found using the identity, 
$\pi_1(\e) = 1 - \pi_2(\e) - \pi_3(\e), \e \in (0, \e_0]$,  and the operational rules for asymptotic expansions  given in Lemma 3$_*$. This yields the asymptotic expansion,   $\pi_1(\e) =   \frac{1}{7} + \frac{6}{147}\e + o_1(\e)$.

Alternatively,  the exclusion of states from the phase space $\XX$ in the opposite order, first state $3$ and then state $2$ or $1$,  and the use of the algorithms described in Lemma  8$_*$ and Theorems 2$_*$ -- 4$_*$ yield the asymptotic expansions for the expected return times $_{\langle 3, 2 \rangle}e_{11}(\e) = E_{11}(\e) = 7 \e + \frac{15}{3} \e^2  +  \ddot{o}_{11}(\e^2)$ and 
$_{\langle 3, 1 \rangle}e_{22}(\e) = E_{22}(\e) = \frac{21}{4} + \frac{15}{4} \e  +  \ddot{o}_{22}(\e)$. Then,  the algorithm described in Theorem 5$_*$ yields the same asymptotic expansions for stationary probabilities  $\pi_1(\e)  = \frac{e_{1}(\e)}{E_{11}(\e)} = \frac{\e + \e^2 + \dot{o}_1(\e^2)}{7\e +  \frac{15}{3} \e^2  + \ddot{o}_{11}(\e^2)} = \frac{1}{7} + \frac{6}{147}\e + o_1(\e)$ and $\pi_2(\e)  = \frac{e_{2}(\e)}{E_{22}(\e)} = \frac{3 + \e + 2 \e^2 + \dot{o}_2(\e^2)}{ \frac{21}{4} + \frac{15}{4} \e  +  \ddot{o}_{22}(\e)} = \frac{4}{7} - \frac{32}{147}\e + o_2(\e)$. 

Note that the Laurent asymptotic expansion for the expectation $E_{22}(\e)$ and, in sequel,  the Taylor asymptotic expansion for the stationary probability $\pi_2(\e)$ are invariant with respect to the choice the sequence of states ${\langle 1, 3 \rangle}$ or ${\langle 3, 1 \rangle}$ for sequential exclusion from the phase space $\XX$. This is consistent with the corresponding invariance statements formulated in Theorems 4$_*$ and 5$_*$.

The coefficients of  the asymptotic expansions $\pi_i(\e) = c_i[0] + c_i[1]\e  + o_i(\e), i = 1, 2, 3$ given above satisfy relations, $c_1[0]  + c_2[0]  + c_3[0] = 1$ and $c_1[1]  + c_2[1]  + c_3[1] = 0$. This is  consistent with the corresponding statement  in Theorem 5$_*$. 

In the example presented above, we did not trace the explicit formulas for remainders $o_1(\e),  o_2(\e)$ and  $o_3(\e)$. However, according to the corresponding statement in  Theorem 
5$_*$,  these  remainders are connected by the following identity $o_1(\e) +  o_2(\e) + o_3(\e) \equiv 0$. 

We would like also to explain  an unexpected,  in some sense, asymptotic behavior of stationary probabilities $\pi_i(\e)$,  in the above example. As a matter of fact,  states $1$ and $2$ are asymptotically absorbing states with non-absorption probabilities of different order, respectively, $O(\e^2)$ and $O(\e)$. While,  state $3$ is a transient  asymptotically non-absorbing  state. This, seems,  should cause convergence of the stationary probability $\pi_1(\e)$ to $1$ and the stationary probabilities $\pi_2(\e)$ and $\pi_3(\e)$ to $0$ as $\e \to 0$, with different rates of convergence. This, however, does not take  place,  and all three probabilities converge to non-zero limits. This is because of the expected  sojourn times $e_1(\e), e_2(\e)$ and  $e_3(\e)$  have orders, respectively, $O(\e), O(1)$ and $O(\e^{-1})$. These expectations compensate  absorption effects for states $1, 2$ and $3$. 

In the above example, computations of explicit upper bounds for remainders in the asymptotic expansions  for stationary probabilities $\pi_1(\e), \pi_2(\e)$ and $\pi_3(\e)$ can also be realized  in the case, where conditions ${\bf D'}$ and ${\bf E'}$ hold instead of conditions  ${\bf D}$ and ${\bf E}$.  We, however, omit this presentation, in order to escape overloading the paper by technical numerical computations.


\begin{thebibliography}{99}
\footnotesize

\bibitem{AFHo} 
{\sc Avrachenkov, K.~E., Filar, J.~A. and Howlett, P.~G.}  (2013). {\em Analytic Perturbation Theory and Its Applications}. SIAM, Philadelphia, PA, xii+372 pp.
\vspace{-2mm} 

\bibitem{BLM} 
{\sc Bini, D.~A., Latouche, G. and Meini, B.} (2005). {\em Numerical Methods for Structured Markov Chains}. Numerical Mathematics and Scientific Computation, Oxford Science Publications,  Oxford University Press, New York,  xii+327 pp.
\vspace{-2mm} 

\bibitem{Cou2}
{\sc Courtois, P.~J.}   (1977). {\em Decomposability. Queueing and Computer System  Applications}. ACM Monograph Series, Academic Press, New
York, xiii+201 pp.
\vspace{-2mm} 

\bibitem{GS3}  
{\sc Gyllenberg, M. and Silvestrov, D.~S.}  (2008). {\em  Quasi-Stationary Phenomena in Nonlinearly Perturbed Stochastic Systems}. De Gruyter 
Expositions in Mathematics, {\bf 44}, Walter de Gruyter, Berlin, ix+579 pp.
\vspace{-2mm} 

\bibitem{Ka10}  
{\sc Kartashov, M.~V.} (1996). {\em Strong Stable Markov Chains}. VSP,
Utrecht and TBiMC, Kiev, 138 pp.
\vspace{-2mm} 

\bibitem{KGM}
{\sc Konstantinov, M.,  Gu, D.~W.,  Mehrmann, V. and Petkov, P.}  (2003). {\em Perturbation Theory for Matrix Equations}.
Studies in Computational Mathematics, {\bf 9}, North-Holland, Amsterdam, xii+429 pp.
\vspace{-2mm} 

\bibitem{Ko2}  
{\sc Korolyuk, V.~S. and Korolyuk, V.~V.} (1999).  {\em Stochastic Models of Systems}. Mathematics and its Applications, 
{\bf 469}, Kluwer, Dordrecht,  xii+185 pp.
\vspace{-2mm} 

\bibitem{KoLi0}
{\sc Koroliuk, V.~S. and Limnios, N.} (2005). {\em Stochastic Systems in
Merging Phase Space}. World Scientific, Singapore, xv+331 pp.
\vspace{-2mm} 

\bibitem{KT3}  
{\sc Korolyuk, V.~S. and Turbin, A.~F.}   (1976). {\em Semi-Markov Processes and its Applications}. Naukova Dumka, Kiev, 184 pp.
\vspace{-2mm} 

\bibitem{KT4}  
{\sc Korolyuk, V.~S. and Turbin, A.~F.}   (1978). {\em Mathematical Foundations of the State Lumping of Large Systems}. Naukova Dumka, Kiev, 218 pp. (English
edition: Mathematics and its Applications, {\bf 264}, Kluwer, Dordrecht, 1993, x+278 pp.).
\vspace{-2mm} 

\bibitem{Sen6}
{\sc Seneta, E.} (2006). {\em  Non-Negative Matrices and Markov Chains}. Springer 
Series in Statistics. Springer, New York,  xvi+287 pp. (A revised reprint of the second (1981) edition).
\vspace{-2mm} 

\bibitem{SS1}  
{\sc Silvestrov, D. and Silvestrov, S.}  (2015). Asymptotic expansions for stationary distributions of perturbed semi-Markov processes. Research Report 2015-9, Department of Mathematics, Stockholm University, 75 pp. and arXiv:1603.03891. 
\vspace{-2mm} 

\bibitem{SS2}  
{\sc Silvestrov,  D. and Silvestrov, S.}  (2016). Asymptotic expansions for stationary distributions of 
nonlinearly perturbed semi-Markov processes. I.  
\vspace{-2mm} 

\bibitem{St9}  
{\sc Stewart, G.~W.}  (1998). {\em Matrix Algorithms. Vol. I. Basic Decompositions}. SIAM,  Philadelphia, PA, xx+458 pp.
\vspace{-2mm} 

\bibitem{St10}  
{\sc Stewart, G.~W.} (2001).  {\em Matrix Algorithms. Vol. II. Eigensystems}. SIAM, Philadelphia, PA, xx+469 pp.
\vspace{-2mm} 

\bibitem{SS1}  
{\sc Stewart, G.~W. and Sun, J.~G.} (1990). {\em Matrix Perturbation Theory}. Computer Science and Scientific Computing. Academic Press,
Boston, xvi+365 pp.
\vspace{-2mm} 

\bibitem{YZ2}  
{\sc Yin, G.~G. and Zhang, Q.}  (2005). {\em Discrete-Time Markov Chains.
Two-Time-Scale Methods and Applications}. Stochastic Modelling and Applied Probability, {\bf 55}, Springer, New York, xix+348 pp.
\vspace{-2mm} 

\bibitem{YZ3}  
{\sc Yin, G.~G. and Zhang, Q.} (2013).  {\em Continuous-Time Markov Chains and Applications. A Two-Time-Scale Approach}. Second edition,  Stochastic Modelling and Applied Probability, {\bf 37}, Springer, New York,  xxii+427 pp. (An extended variant of the first (1998)  edition).
\vspace{-2mm} 




\end{thebibliography}
\end{document}